\documentclass[12pt]{amsart}
\usepackage[T1]{fontenc}
\usepackage[english]{babel}
\usepackage{mathtools}
\usepackage{amstext, amssymb}
\usepackage{mathrsfs}
\usepackage{nicematrix}
\usepackage{slashed}
\usepackage{accents}
\usepackage{enumerate}
\usepackage{xfrac}
\usepackage[all]{xy}

\theoremstyle{plain}
\newtheorem{theorem}{Theorem}[section]
\newtheorem{proposition}[theorem]{Proposition}
\newtheorem{corollary}[theorem]{Corollary}
\newtheorem{lemma}[theorem]{Lemma}

\theoremstyle{definition}
\newtheorem{definition}[theorem]{Definition}
\newtheorem{remark}[theorem]{Remark}

\theoremstyle{remark}
\newtheorem{notation}[theorem]{Notation}

\numberwithin{equation}{section}

\DeclareMathOperator{\Hom}{\mathrm{Hom}}
\DeclareMathOperator{\End}{\mathrm{End}}
\DeclareMathOperator{\dom}{\mathrm{Dom}}
\DeclareMathOperator{\sgn}{\mathrm{sgn}}
\DeclareMathOperator{\Str}{\mathrm{Str}}

\newcommand{\IMM}{\mathscr{M}}
\newcommand{\ILL}{\mathscr{L}}
\newcommand{\IHH}{\mathscr{H}}
\newcommand{\IFF}{\mathscr{F}}

\newcommand{\IN}{\mathbb{N}}
\newcommand{\IZ}{\mathbb{Z}}
\newcommand{\IP}{\mathbb{P}}
\newcommand{\IC}{\mathbb{C}}
\newcommand{\IR}{\mathbb{R}}

\newcommand{\Id}{\mathrm{d} }

\newcommand{\pa}{\slash \hspace{-1mm}\slash}
\newcommand{\nn}{\nonumber}

\def\transport{/\hspace*{-3pt}/}

\def\Xint#1{\mathchoice {\XXint\displaystyle\textstyle{#1}}%
{\XXint\textstyle\scriptstyle{#1}}%
{\XXint\scriptstyle\scriptscriptstyle{#1}}%
{\XXint\scriptscriptstyle\scriptscriptstyle{#1}}%
\!\int}
\def\XXint#1#2#3{{\setbox0=\hbox{$#1{#2#3}{\int}$} \vcenter{\hbox{$#2#3$}}\kern-.5\wd0}}
\def\dashint{\Xint-}

\allowdisplaybreaks

\title[]{Fermionic Dyson expansions and stochastic Duistermaat-Heckman localization on loop spaces}

\author{Batu Güneysu}
\author{Jonas Miehe}

\begin{document}
	\begin{abstract}
		Given a self-adjoint operator $H\geq 0$ and (appropriate) densely defined and closed operators $P_{1},\dots, P_{n}$ in
		a Hilbert space $\IHH$, we provide a systematic study of bounded operators given by iterated integrals
		\begin{align}\label{oh}
                \int_{\{ 0\leq s_1\leq \dots\leq s_n\leq t\}}\mathrm{e}^{-s_1H}P_{1}\mathrm{e}^{-(s_2-s_1)H}P_{2}\cdots \mathrm{e}^{-(s_n-s_{n-1})H}P_{n} \mathrm{e}^{-(t-s_n)H}\, \Id s_{1} \ldots \Id s_{n},\quad t>0.
		\end{align}
		These operators arise naturally in noncommutative geometry and the geometry of loop spaces. Using Fermionic calculus,
		we give a natural construction of an enlarged Hilbert space $\IHH^{(n)}$ and an analytic semigroup
		$\mathrm{e}^{-t (H^{(n)}+P^{(n)} )}$ thereon, such that $\mathrm{e}^{-t (H^{(n)}+P^{(n)} )}$ composed from the left with
		(essentially) a Fermionic integration gives precisely the iterated operator integral (\ref{oh}). This formula allows
		to establish important regularity results for the latter, and to derive a stochastic representation for it, in case
		$H$ is a covariant Laplacian and the $P_{j}$'s are first-order differential operators. Finally, with $H$ given as
		the square of the Dirac operator on a spin manifold, this representation is used to derive a stochastic refinement of
		the Duistermaat-Heckman localization formula on the loop space of a spin manifold.
	\end{abstract}

	\maketitle

	\section{Summary of the main results}

	Let $\IHH$ be a complex Hilbert space, $H\geq 0$ be a self-adjoint operator in $\IHH$, and let $P_{1},\dots,P_{n}$
	be densely defined and closed operators in $\IHH$. In this paper, we are going to systematically study bounded operators in $\IHH$ of the form
	\begin{align*}
		 & \Phi^{H}_{t}(P_{1},\dots,P_{n}) \\
		 & \coloneqq\int_{\{ 0\leq s_1\leq \dots\leq s_n\leq t\}}\mathrm{e}^{-s_1H}P_{1}\mathrm{e}^{-(s_2-s_1)H}P_{2}\cdots \mathrm{e}^{-(s_n-s_{n-1})H}P_{n} \mathrm{e}^{-(t-s_n)H}\, \Id s_{1} \ldots \Id s_{n},\quad t>0.
	\end{align*}
	While, in view of the Dyson expansion
	\begin{equation*}
		\mathrm{e}^{-t(H+P)}=\mathrm{e}^{-tH}+\sum^{n}_{l=1}(-1)^{l}\Phi^{H}_{t}\underbrace{(P,\dots,P)}_{\text{$l$ times}},
	\end{equation*}
	the case $P_{1}=\dots = P_{n}=P$ appears naturally in perturbation theory, the general case plays a fundamental role in
	noncommutative geometry.

	\subsection{Differential Graded JLO Cocycle}
	The theory of even $\vartheta$-summable Fredholm modules over locally convex unital algebras has been initiated by Connes
	\cite{connes}, Jaffe / Lesniewski / Osterwalder \cite{jlo} and Getzler / Szenes \cite{gs} around 1990. A refinement of
	this theory to \emph{differential graded} algebras (DGAs) has been constructed recently in \cite{gl} and generalized in
	\cite{miehe}. Here, starting from a locally convex unital DGA $\Omega$, an \emph{even $\vartheta$-summable Fredholm
	module $\IMM$ over $\Omega$} is given by
	\begin{itemize}
		\item a $\IZ_{2}$-graded Hilbert space $\IHH$,

		\item an odd unbounded self-adjoint operator $D$ in $\IHH$ (an abstract Dirac operator) with a trace-class heat semigroup
			$\mathrm{e}^{-tD^2}$ for all $t>0$,

		\item an even and continuous linear map (an abstract Clifford multiplication)
			\begin{equation*}
				\mathbf{c}\colon\Omega\longrightarrow\ILL(\IHH),
			\end{equation*}
			which, in addition to some analytic conditions, satisfies
			\begin{align*}
				\mathbf{c}(1_{\Omega})=\mathrm{id}_{\IHH},\quad [D,\mathbf{c}(f)]=\mathbf{c}(\Id f),\quad \mathbf{c}(f\omega)=\mathbf{c}(f)\mathbf{c}(\omega),\quad \mathbf{c}(\omega f)=\mathbf{c}(\omega)\mathbf{c}(f),\quad f\in \Omega^{0}, \omega\in \Omega.
			\end{align*}
	\end{itemize}
	In particular, $\mathbf{c}$ is only a representation of the degree zero part $\Omega^{0}$ of $\Omega$, and the $\IZ_{2}$-graded
	commutators $[D, \mathbf{c}(\omega)]$ are not assumed to be bounded in case $\deg(\omega)>0$, which leads to certain
	difficulties in comparison to the ungraded theory.

	One then considers the \emph{acyclic extension} $\Omega_{\mathbb{T}}$ of $\Omega$,
	with $\Omega_{\mathbb{T}}$ defined as the locally convex unital DGA of all $\omega=\omega'+u\omega''$ with
	$\omega',\omega''\in\Omega$ and $u$ a formal variable of degree $-1$ with $u^{2}=0$. With the \emph{cyclic complex}
	\begin{equation*}
		\mathsf{C}^{*}(\Omega_{\mathbb{T}})\coloneqq\bigoplus^{\infty}_{n=0}\Omega_{\mathbb{T}}\otimes (\Omega_{\mathbb{T}}/\IC)^{\otimes
		n},
	\end{equation*}
	the Chern character of $\IMM$ -- which in view of \cite{jlo} can be regarded as a differential graded JLO cocycle -- is
	a cocycle
	\begin{equation*}
		\mathrm{Ch}_{\mathbb{T}}(\IMM)\in \mathsf{C}_{*}(\Omega_{\mathbb{T}})=\Hom(\mathsf{C}^{*}(\Omega_{\mathbb{T}}),\IC),
	\end{equation*}
	which is defined as follows: firstly, given $\omega_{1}, \ldots, \omega_{k} \in \Omega_{\mathbb{T}}$, one gets the densely
	defined closed operators
	\begin{align*}
		\begin{split}P(\omega_{1})&\coloneqq [D, \mathbf{c}(\omega_{1}')] - \mathbf{c}(\Id \omega_{1}') + \mathbf{c}(\omega_{1}''), \\ P(\omega_{1}, \omega_{2})&\coloneqq (-1)^{\deg(\omega_1')}\big(\mathbf{c}(\omega_{1}' \omega_{2}') - \mathbf{c}(\omega_{1}') \mathbf{c}(\omega_{2}')\big), \\ P(\omega_{1}, \ldots, \omega_{k})&\coloneqq 0, \quad k \geq 3.\end{split}
	\end{align*}
	Then, with $\mathcal{P}_{m, n}$ the set of ordered partitions of length $m$ of $\{1, \ldots, n\}$, one can define
	\begin{align}\label{Definition Chern 2}
            \mathrm{Ch}_{\mathbb{T}}(\IMM)(\omega_{0}, \ldots, \omega_{n}) & \coloneqq \Str\Big(\mathbf{c}(\omega_{0}')\sum_{m = 1}^{n} (-1)^{m} \!\!\!\sum_{(I_1, \ldots, I_m) \in \mathcal{P}_{m,n}}\Phi_{t=1}^{D^2}\big(P(\omega_{I_1}), \ldots, P(\omega_{I_m})\big)\Big),
	\end{align}
	where
	\begin{equation*}
		\text{$\omega_{I_a}\coloneqq(\omega_{i+1},\dots, \omega_{i+i'})$, if $I_{a}= \{j \colon i < j \leq i + i'\}$ for some $i$,
		$i'$}.
	\end{equation*}
	A fundamental result is the homotopy invariance \cite{gl, miehe} of the differential graded JLO cocycle, which implies
	\begin{align}
		\label{homo}\lim_{t\to 0+}\mathrm{Ch}_{\mathbb{T}}(\IMM_{t})(\alpha)=\mathrm{Ch}_{\mathbb{T}}(\IMM)(\alpha)\quad \text{for every cycle $\alpha\in \mathsf{C}^{*}(\Omega_{\mathbb{T}}(X))$}.
	\end{align}
	It is also of fundamental importance that $\mathrm{Ch}_{\mathbb{T}}(\IMM_{t})$ satisfies an analytic growth condition in the
	spirit of \cite{connes,jlo}, so that one can actually pair this linear form with certain infinite chains obtained by completing
	$\mathsf{C}^{*}(\Omega_{\mathbb{T}}(X))$ in a natural way.

	As will be explained in Section \ref{wyaq} below, apart from being merely a refinement of the ungraded theory, the graded
	theory turns out to be essential in the context of loop space geometry.

	\subsection{$\Phi^{H}_{t}(P_{1},\dots,P_{n})$ and Fermionic Dyson Expansions}
	Returning to the general situation, note first that it is a simple consequence of spectral calculus that under the
	assumption
	\begin{align*}
		P_{j}(H+1)^{-a_j}\in \ILL(\IHH)\quad\text{for some $a_{j}\in (0,1)$,}
	\end{align*}
	the operator $\Phi^{H}_{t}(P_{1},\dots,P_{n})$ can be defined as a Bochner integral in the Banach space $\ILL(\IHH)$ (cf.
	Lemma 2.3). The central result of this paper deals the following question:
        \begin{center}
            \emph{In what sense is $\Phi^{H}_{t}(P_{1},\dots,P_{n})$ algebraically induced by an operator semigroup?}
        \end{center}
	Our main result in this regard states that one can use algebraic methods inspired by quantum field
	theory (Fermionic calculus) to canonically construct such a (analytic) operator semigroup, which when composed with a variant
	of the \emph{Fermionic} or \emph{Berezin integral} reproduces $\Phi^{H}_{t}(P_{1},\dots,P_{n})$. Concretely, let
	\begin{equation*}
		\IHH^{(n)}\coloneqq (\Lambda_{n}\otimes \IHH)^{n}=\bigoplus^{n}_{j=1}(\Lambda_{n}\otimes \IHH).
	\end{equation*}
	Then one obtains the Fermionic integral (cf. Section \ref{ddqay})
	\begin{equation*}
		\dashint_{ \Lambda_n^* \otimes \ILL(\IHH^n)}\colon \ILL(\IHH^{(n)})\longrightarrow \Lambda_{n}^{*}\otimes \ILL(\IHH^{n}),
	\end{equation*}
	as well as the linear map
	\begin{equation*}
		\pi^{(n)}_{\IHH}\colon \Lambda_{n}^{*}\otimes \ILL(\IHH^{n})\longrightarrow \ILL( \IHH), \quad \alpha \otimes (a_{ij})_{i,j=1,\dots,n}
		\longmapsto a_{nn}.
	\end{equation*}

	By letting $H$ act diagonally, we obtain a selfadjoint operator $H^{(n)}\geq 0$ in $\IHH^{(n)}$. The action of $P$ is more
	involved, and is ultimately given by the closed operator
	\begin{equation*}
		P^{(n)}\coloneqq
		\begin{pmatrix}
                & & & \hat{\theta}_{n} \otimes P_{n} \\
			\hat{\theta}_{n-1}\otimes P_{n-1} & & & \\
			  & \ddots & & \\
			  & & \hat{\theta}_{1}\otimes P_{1} &
		\end{pmatrix},
	\end{equation*}
	where for every $\alpha \in \Lambda_{n}$ the symbol $\hat{\alpha}$ denotes the induced multiplication operator in
	$\Lambda_{n}$. With these preparations, Theorem \ref{main}, states that $H^{(n)}+P^{(n)}$ generates an analytic semigroup
	in $\IHH^{(n)}$, to be denoted by $\mathrm{e}^{-t(H^{(n)}+P^{(n)})}$, which is given by the (finite) fermionic Dyson expansion
	\begin{equation*}
		\mathrm{e}^{-t(H^{(n)}+P^{(n)})}=\mathrm{e}^{-tH^{(n)}}+\sum^{n}_{l=1}(-1)^{l}\Phi^{H^{(n)}}_{t}\underbrace{(P^{(n)},\dots,P^{(n)})}
		_{\text{$l$ times}},
	\end{equation*}
	and moreover
	\begin{align}
		\label{mai}\Phi_{t}^{H}(P_{1},\dots, P_{n})=(-1)^{n} \pi^{(n)}_{\IHH}\dashint_{ \Lambda_n^* \otimes \ILL(\IHH^n)}\mathrm{e}^{-t(H^{(n)}+P^{(n)})}.
	\end{align}
	An immediate consequence of the latter formula is that the map $t\mapsto \Phi^{H}_{t}(P_{1},\dots,P_{n})$ is norm smooth
	on $(0,\infty)$ and extends strongly continuously to $0$ at $t=0$. Moreover, one obtains the smoothing property
	\begin{equation*}
		\Phi^{H}_{t}(P_{1},\dots,P_{n})\colon\IHH\longrightarrow \dom(H),
	\end{equation*}
        while in general $\Phi^{H}_{t}(P_{1},\dots,P_{n})$ does not even map
	\begin{equation*}
		\bigcap_{r=1}^{\infty} \dom(H^{r})\longrightarrow \dom(H^{2}).
	\end{equation*}

	\subsection{$\Phi^{H}_{t}(P_{1},\dots,P_{n})$ for Geometric Operators}

	We then move on to focus on operators that arise in geometry. Here, we consider a covariant Schrödinger operator
	$H=\nabla^{\dagger} \nabla/2+W$, with $\nabla$ a metric connection and $W$ self-adjoint of zeroth-order, both acting
	on a metric vector bundle $E\to X$ over a closed Riemannian manifold $X$. The $P_{j}$'s are assumed to be differential
	operators of order $\leq 1$. In this situation, we obtain e. g. the existence of a jointly smooth heat kernel
	\begin{equation*}
		(0,\infty)\times X\times X \ni (t,x,y)\longmapsto \Phi^{\nabla^\dagger \nabla/2+W}_{t}(P_{1},\dots,P_{n})(x,y)\in \Hom
		(E_{y},E_{x})\subset E\boxtimes E^{*},
	\end{equation*}
	which admits a strong asymptotic expansion as $t\to 0+$ (cf. Theorem \ref{expan}). A fundamentally new consequence of
	(\ref{mai}) in the geometric setting is given by the following probabilistic representation formula, which is due to Theorem \ref{Feynman-Kac},
	\begin{align}
		\label{ooh} & \Phi^{\nabla^\dagger\nabla/2 +W}_{t}(P_{1},\dots,P_{n})(x,y) \\
		\nn & \quad = p(t,x,y)\; \mathbb{E}\Big[ \mathscr{W}_{\nabla}^{x,y;t}(t) \; \int_{\{ 0\leq s_1\leq \dots\leq s_n\leq t\}}\delta\Psi^{x,y;t}_{\nabla,P_1}(s_{1})\cdots \delta\Psi^{x,y;t}_{\nabla,P_n}(s_{n}) \;\pa^{x,y;t}_{\nabla}(t)^{-1}\Big],
	\end{align}
	where
	\begin{itemize}
		\item[i)] $p(t,x,y)$ denotes the scalar heat kernel,

		\item[ii)] $\mathscr{W}_{\nabla}^{x,y;t}(\bullet)$ is the usual multiplicative functional induced by the potential $W$,
			calculated with respect to $\nabla$ and any adapted Brownian bridge $\mathsf{B}_{\bullet}^{x,y;t}$ starting in $x$
			and ending in $y$ at the time $t$,

		\item[iii)] the process
			\begin{equation*}
				\int_{\{ 0\leq s_1\leq \dots\leq s_n\leq \bullet\}}\delta\Psi^{x,y;t}_{\nabla,P_1}(s_{1})\cdots \delta\Psi^{x,y;t}
				_{\nabla,P_n}(s_{n})
			\end{equation*}
			is an iterated It\^{o} line integral along $\mathsf{B}_{\bullet}^{x,y;t}$, which is calculated using the
			decomposition with respect to $\nabla$ of each $P_{j}$ into a pure first-order part and a pure zeroth-order part,

		\item[iv)] $\pa^{x,y;t}_{\nabla}(\bullet)^{-1}$ denotes the inverse stochastic parallel transport, calculated with
			respect to $\nabla$ and $\mathsf{B}_{\bullet}^{x,y;t}$.
	\end{itemize}

	\subsection{$\Phi^{H}_{t}(P_{1},\dots,P_{n})$ and the Geometry of Loop Spaces}
	\label{wyaq}

	Formula (\ref{ooh}) plays a fundamental role in the geometry of the loop space $LX$ over a closed even-dimensional
	Riemannian spin manifold $X$. To explain this connection, consider the even $\vartheta$-summable Fredholm module $\IMM^{X}$
	over $\Omega(X)$, given by the Dirac operator $\slashed{D}:=2^{-\frac{1}{2}}D$ acting on the Hilbert space of $L^{2}$-spinors on $X$, with the Clifford multiplication $\slashed{\mathbf{c}}(\omega):=2^{-\frac{\mathrm{deg}(\omega)}{2}}\mathbf{c}$, where $D$ resp. $\mathbf{c}$ denote the usual Dirac operator resp. Clifford multiplication. As recorded in \cite{h1,h2,gl,bcg}, for every $\alpha\in \mathsf{C}^{*}(\Omega_{\mathbb{T}}
	(X))$ the number $\mathrm{Ch}_{\mathbb{T}}(\IMM^{X})(\alpha)\in\IC$ can be interpreted as a definition of the Gaussian
	integral of the Chen differential form on $LX$ induced by $\alpha$. This definition is functional analytic in its
	spirit, much in the way one can define the Wiener measure without referring to paths at all (see e. g. Section \ref{appe}).
	This functional analytic construction satisfies all desired properties of an integration theory for differential forms
	on $LX$. Most notably, keeping in mind the natural $S^{1}$-action on $LX$ given by rotating loops, there holds a
	version of the Duistermaat-Heckman localization formula \cite{gl,lude2}, which states that
	\begin{align}
		\label{uzyx}\mathrm{Ch}_{\mathbb{T}}(\IMM^{X})(\alpha)= \int_{X}\hat{A}(X)\wedge h(\alpha)\quad\text{for every cycle $\alpha\in \mathsf{C}^{*}(\Omega_{\mathbb{T}}(X))$.}
	\end{align}
	Here,
	\begin{itemize}
		\item $\hat{A}(X)$ denotes the Chern-Weil representative of the $\hat{A}$-genus,

		\item the linear map $h$ is defined via
			\begin{equation*}
				h\colon\mathsf{C}^{*}(\Omega_{\mathbb{T}}(X))\longrightarrow \Omega(X), \quad h(\omega_{0},\dots,\omega_{n})\coloneqq\frac{(-1)^{n}
				\, 2^{2n}}{n! (2 \pi \sqrt{-1})^{\frac{\dim(X)}{2}}}\;\omega'_{0}\wedge\omega''_{1}\wedge\cdots \wedge \omega_{n}
				''.
			\end{equation*}
	\end{itemize}
	While to a large degree this formula -- several decades after the heuristic considerations by Atiyah \cite{atiyah}, 
        Bismut \cite{bismut2}, and Witten \cite{witten} -- justifies mathematically the loop space localization proof of the  Atiyah-Singer index theorem, a major open question, that could only be answered in special cases so far \cite{h2,boldt}, is:
        \begin{center}
            \emph{Is there a formula for $\mathrm{Ch}_{\mathbb{T}}(\IMM^{X})$ which is given in terms of an actual integration on
            a space of loops, and which makes the effects leading to the localization formula transparent at the level of loops?}
        \end{center}
	
	Formula (\ref{ooh}) leads to a complete answer to this question at the level of Brownian loops: consider the homotopy $\IMM
	^{X}_{t}$ of even $\vartheta$-summable Fredholm modules which is given for $t\in [0,1]$ by replacing $\slashed{D}$ with $t^{\frac{1}{2}}
	\slashed{D}$ and $\slashed{\mathbf{c}}(\omega)$ with $ t^{\frac{\mathrm{deg}(\omega)}{2}}\slashed{\mathbf{c}}(\omega)$. In view of
	the Lichnerowicz formula
	\begin{equation*}
		\slashed{D}^2 = \nabla^{\dagger}\nabla/2 + \mathrm{scal}/8,
	\end{equation*}
        and with $\mu$ the Riemannian volume measure, we get the fundamental relation\footnote{Note that we use the $\vartheta$-summable Fredholm module $\IMM^X$ induced by a rescaled Dirac operator and a rescaled Clifford multiplication, in order to be consistent with the conventions in stochastic analysis. This is the reason for the powers of two appearing in \eqref{erad} compared to \eqref{Definition Chern 2}.}
	\begin{align}\nn
            & \mathrm{Ch}_{\mathbb{T}}(\IMM^{X}_{t})(\omega_{0}, \ldots, \omega_{n}) \\
            \label{erad}  & = \int_{X} \mathbb{E}\Big[(t/2)^{ -\frac{n}{2}+\frac{1}{2}\sum^n_{j=0}\mathrm{deg}(\omega'_j) }\Str_{x} \Big(\sum_{m = 1}^{n} (-2)^{m} \sum_{(I_1,
            \ldots, I_m) \in \mathcal{P}_{m,n}}\mathbf{c}_{x}(\omega_{0}')p(t,x,x) \times \\
            \nn & \quad \times \mathrm{e}^{-\frac{1}{8}\int_0^t \mathrm{scal}(\mathsf{B}_s^{x,x;t}) \, \Id s}\int_{t\sigma_n}\delta\Psi^{x,x;t}_{\nabla,
            P(\omega_{I_1})}(s_{1})\cdots \delta\Psi^{x,x;t}_{\nabla, P(\omega_{I_m})}(s_{m}) \pa^{x,x;t}_{\nabla}(t)^{-1}\Big)\Big] \, \Id\mu(x).
	\end{align}
	Our main result in this context, Theorem \ref{local}, states the following stochastic localization formula: uniformly
	in $x\in X$, the limit
	\begin{align}\label{limi}
            F_{x}(\omega_{0},\dots,\omega_{n}) \coloneqq &\>\lim_{t\to 0+}(t/2)^{ -\frac{n}{2}+\frac{1}{2}\sum^n_{j=0}\mathrm{deg}(\omega'_j) }\Str_{x} \Big(\sum_{m = 1}^{n} (-2)^{m} \sum_{(I_1,\ldots, I_m) \in \mathcal{P}_{m,n}}\mathbf{c}_{x}(\omega_{0}')p(t,x,x) \times \\
            \nn & \quad \times \mathrm{e}^{-\frac{1}{8}\int_0^t \mathrm{scal}(\mathsf{B}_s^{x,x;t}) \, \Id s}\int_{t\sigma_n}\delta\Psi^{x,x;t}_{\nabla,
            P(\omega_{I_1})}(s_{1})\cdots \delta\Psi^{x,x;t}_{\nabla, P(\omega_{I_m})}(s_{m}) \pa^{x,x;t}_{\nabla}(t)^{-1}\Big)
	\end{align}
 
	exists in $L^{b}(\IP)$ for all $b\in [1,\infty)$. Moreover, with $\underline{\mu}$ the Riemannian volume form on $X$,
	the complex random field $F_{\bullet}(\omega_{0},\dots,\omega_{n})$ satisfies
	\begin{equation*}
		\mathbb{E}\left[F_{\bullet}(\omega_{0},\dots,\omega_{n})\right]\cdot \underline{\mu}= \frac{(-1)^{n} \, 2^{2n}}{n! (2
		\pi \sqrt{-1})^{\frac{\dim(X)}{2}}}\hat{A}(X)\wedge \omega'_{0}\wedge\omega''_{1}\wedge\cdots \wedge \omega_{n}'' \quad\text{as volume forms}.
	\end{equation*}
	In view of the homotopy invariance (\ref{homo}) and formula (\ref{erad}), these results can be seen as a (localized) stochastic
	refinement of (\ref{uzyx}), making a bridge between noncommutative geometry, probability theory and the geometry of loop
	spaces.

	Finally, concerning the effects leading to the existence of the limit (\ref{limi}), our proof reveals that the following two facts ultimately combine to a stochastic localization:
	\begin{itemize}
		\item \emph{Analytic cancellations}: let $F_{x;t}(\omega_{0},\dots,\omega_{n})$ denote the complex random variable given
			by the right hand side (\ref{limi}) before taking the limit. Then the on-diagonal singularity of the scalar heat
			kernel, Hsu's \cite{hsu} asymptotic Taylor expansion of the Brownian holonomy $\pa^{x,x;t}_{\nabla}(t)^{-1}$ (cf. Proposition \ref{Limit
			parallel transport}), and careful $L^{b}(\IP)$ estimates (cf. Lemma \ref{It\^{o}}) of the iterated It\^{o}
			integrals show that certain summands in $F_{x;t}(\omega_{0},\dots,\omega_{n})$ must vanish as $t\to 0+$ after
			taking expectation.

		\item \emph{Combinatoric cancellations}: using the Patodi filtration of the endomorphisms of the spinor bundle (cf. Definition
			\ref{patodifilt}), one obtains that certain summands in $F_{x;t}(\omega_{0},\dots,\omega_{n})$ must vanish already
			at the level of based Brownian loops, that is, before taking expectation and the limit $t\to 0+$. This reasoning
			relies on rather subtle calculations involving the iterated It\^{o} integrals and constitutes the hard part of the
			proof.
	\end{itemize}
 
\vspace{2mm}
        \textbf{Acknowledgements} The authors would like to thank Sebastian Boldt and Elton Hsu for several helpful discussions.
\vspace{2mm}

\vspace{2mm}
        \emph{To our beloved parents Mihriban \& Erdolon Güneysu, Verena \& Michael Miehe.         }
\vspace{2mm}

	\section{Basic Properties of $\Phi^{H}_{t}(P_{1},\dots,P_{n})$}
	\label{ddqay}

	Given $n\in\IN$, let $\Lambda_{n}$ denote the exterior algebra of $\IC^{n}$ and let $\theta_{1},\dots,\theta_{n}$ be the
	standard basis of $\IC^{n}$. Then one has $\theta_{j}^{2}=0$, $\theta_{i}\theta_{j}+\theta_{j}\theta_{i} = 0$ for $i\neq
	j$ in $\Lambda_{n}$. Furthermore, the products
	\begin{align}
		\label{onb}\theta_{j_1}\cdots \theta_{j_k}, \quad \text{where}\quad 1\leq j_{1}<\ldots <j_{k}\leq n,\quad 1\leq k\leq n,
	\end{align}
	form a basis of $\Lambda_{n}$ so that every $\alpha \in \Lambda_{n}$ can be uniquely written as
	\begin{equation*}
		\alpha=\sum^{n}_{k=0}\sum_{1\leq j_1<\ldots <j_k\leq n}\alpha_{j_1,\ldots,j_k}\theta_{j_1}\cdots \theta_{j_k},\quad \alpha
		_{j_1,\ldots,j_k}\in \IC.
	\end{equation*}
	A scalar product on $\Lambda_{n}$ is given by declaring (\ref{onb}) to be an orthonormal basis.

	It follows that, given a complex linear space $E$, every $\alpha \in \Lambda_{n}\otimes E$ can be uniquely written as
	\begin{equation*}
		\alpha=\sum^{n}_{k=0}\sum_{1\leq j_1<\ldots <j_k\leq n}\theta_{j_1}\cdots \theta_{j_k}\otimes \alpha_{j_1,\ldots,j_k}
		,\quad \alpha_{j_1,\ldots,j_k}\in E.
	\end{equation*}
	\begin{definition}
		The \emph{Fermionic or Berezin integral on $E$} is the linear map
		\begin{equation*}
			\dashint_{ E}\colon \Lambda_{n}\otimes E\longrightarrow E, \quad \alpha \longmapsto \alpha_{1,\ldots,n}.
		\end{equation*}
	\end{definition}

	Let $\IHH$ be a complex Hilbert space and let $\ILL(\IHH)$ denote the space of bounded operators in $\IHH$. We record that
	for a map from an open subset of $\IR$ (resp. $\IC$) to $\ILL(\IHH)$, the notions strongly/norm smooth (resp. analytic)
	are equivalent -- a simple consequence of the uniform boundedness principle.

	Let $H\geq 0$ be a self-adjoint operator in $\IHH$. The semigroup $\mathrm{e}^{-tH}$, $t>0$, is defined via spectral
	calculus and is an analytic semigroup of angle $\pi/2$. As such, it holds that
	\begin{align}
		\mathrm{Ran}(\mathrm{e}^{-t H})\subset C^{\infty}(H)\coloneqq\bigcap^{\infty}_{l=1}\dom(H^{l})\quad\text{for all $t> 0$},
	\end{align}
	Noting that $C^{\infty}(H)$ is a dense subspace of $\IHH$. We refer the reader to \cite{nagel}, section II.4, for
	details on analytic semigroups, and record for future reference the following auxiliary result (cf. Theorem 1.19 in \cite{kato}):

	\begin{lemma}
		\label{aux} Let $g \colon [0,\infty)\to \IHH$ be a $C^{1}$-function. Then the function
		\begin{equation*}
			[0,\infty)\ni t\longmapsto \int^{t}_{0} \mathrm{e}^{-sH}g(t-s) \, \Id s\in \IHH,
		\end{equation*}
		defined via Bochner integration on $\IHH$, is $C^{1}$ with
		\begin{align*}
			\int^{t}_{0} \mathrm{e}^{-sH}g(t-s) \, \Id s\in \dom (H)\quad \text{for all $t\geq 0$,}
		\end{align*}
		and
		\begin{align*}
			(d/dt)\int^{t}_{0} \mathrm{e}^{-sH}g(t-s) \, \Id s & = -H \int^{t}_{0} \mathrm{e}^{-sH}g(t-s) \, \Id s+g(t)\quad \text{for all $t\geq 0$,} \\
			\int^{t}_{0} \mathrm{e}^{-sH}g(t-s) \, \Id s\Bigg\vert_{t=0} & = 0.
		\end{align*}
	\end{lemma}
	
	Assume now $P_{1}, \dots, P_{n}$ are densely defined closed operators in $\IHH$ such that
	\begin{align}\label{ass}
		P_{j}(H+1)^{-a_j}\in \ILL(\IHH)\quad\text{for some $a_{j}\in (0,1)$.}
	\end{align}

	The assumption (\ref{ass}), which is automatically satisfied if $P_{j}$ is bounded, implies that $P_{j}(H+1)^{-1}$ is
	bounded, too, $\dom(H)\subset \dom(P_{j})$, and that $P_{j}\mathrm{e}^{-rH}$ is bounded with
	\begin{align}
		\label{wsx2}\left\|P_{j}\mathrm{e}^{-rH}\right\|=\left\|P_{j}(H+1)^{-a_j}(H+1)^{a_j}\mathrm{e}^{-rH}\right\|\leq \left\|P_{j}(H+1)^{-a_j}\right\|\mathrm{e}^{r}r^{-a_j}\quad\text{for all $r>0$}.
	\end{align}
	It follows that for $t>0$, setting
	\begin{equation*}
		t\sigma_{n}\coloneqq\{s\in\IR^{n} \colon 0\leq s_{1}\leq \dots\leq s_{n}\leq t\},
	\end{equation*}
	the map
	\begin{equation*}
		t\sigma_{n}\ni s\longmapsto \mathrm{e}^{-s_1H}P_{1}\mathrm{e}^{-(s_2-s_1)H}P_{2}\cdots \mathrm{e}^{-(s_n-s_{n-1})H}P_{n}
		\mathrm{e}^{-(t-s_n)H}\in \ILL(\IHH)
	\end{equation*}
	is well-defined almost everywhere. Moreover, from $\left\|\mathrm{e}^{-s_1H}\right\|\leq 1$ and $(\ref{wsx2})$, one obtains
	the following estimate
	\begin{align}\label{Fundamental estimate n-simplex}
            \begin{split}
                &\int_{t\sigma_n}\left\|\mathrm{e}^{-s_1H}P_{1}\mathrm{e}^{-(s_2-s_1)H}P_{2}\cdots \mathrm{e}^{-(s_n-s_{n-1})H}P_{n} \mathrm{e}^{-(t-s_n)H}\right\| \, \Id s_{1}\ldots \Id s_{n}\\&\leq c(P_{1},\dots, P_{n}) \mathrm{e}^{nt}\int_{t\sigma_n}(s_{2}-s_{1})^{-a_1}\cdots (s_{n}-s_{n-1})^{-a_{n-1}}(t-s_{n})^{-a_n}\, \Id s_{1}\ldots \Id s_{n}\\&=c(P_{1},\dots, P_{n}) \mathrm{e}^{nt}t^{n-\sum^n_{j=1}a_j}\int_{\sigma_n}(s_{2}-s_{1})^{-a_1}\cdots (1-s_{n})^{-a_n}\, \Id s_{1}\ldots \Id s_{n}<\infty,
            \end{split}
	\end{align}
	where
	\begin{equation*}
		c(P_{1},\dots, P_{n})\coloneqq\left\|P_{1}(H+1)^{-a_1}\right\|\cdots \left\|P_{n}(H+1)^{-a_n}\right\|.
	\end{equation*}

	We can thus define the Bochner integral
	\begin{equation*}
		\Phi^{H}_{t}(P_{1},\dots,P_{n})\coloneqq\int_{t\sigma_n}\mathrm{e}^{-s_1H}P_{1}\mathrm{e}^{-(s_2-s_1)H}P_{2}\cdots \mathrm{e}
		^{-(s_n-s_{n-1})H}P_{n} \mathrm{e}^{-(t-s_n)H}\, \Id s_{1} \ldots \Id s_{n}\in \ILL(\IHH).
	\end{equation*}

	We clarify that for $n=1$ we understand
	\begin{equation*}
		\Phi^{H}_{t}(P_{1})=\int^{t}_{0}\mathrm{e}^{-sH}P_{1} \mathrm{e}^{-(t-s)H}\, \Id s,
	\end{equation*}
	and it will also we convenient to define
	\begin{align}
		\label{con}\Phi^{H}_{t}(\emptyset)\coloneqq \mathrm{e}^{-tH},
	\end{align}
	as well as
	\begin{equation*}
		\Phi^{H}_{0}(P_{1},\dots,P_{n})\coloneqq 0.
	\end{equation*}
	Then the above estimate can be restated to be:

	\begin{lemma}
		\label{con2} Under (\ref{ass}), for all $t\geq 0$ one has
		\begin{align*}
			\left\|\Phi^{H}_{t}(P_{1},\dots,P_{n})\right\|\leq C(P_{1},\dots, P_{n}) \mathrm{e}^{nt}t^{n-\sum^n_{j=1}a_j},
		\end{align*}
		where
		\begin{equation*}
			C(P_{1},\dots, P_{n})\coloneqq c(P_{1},\dots, P_{n})\int_{\sigma_n}(s_{2}-s_{1})^{-a_1}\cdots (s_{n}-s_{n-1})^{-a_{n-1}}
			(1-s_{n})^{-a_n}\, \Id s_{1}\cdots \Id s_{n}.
		\end{equation*}
	\end{lemma}

	To examine $\Phi^{H}_{t}(P_{1},\dots,P_{n})$ further, let us now define a densely defined closed operator $P^{(n)}$ in
	the Hilbert space
	\begin{equation*}
		\IHH^{(n)}\coloneqq \Lambda_{n}\otimes \IHH\otimes \IC^{n}=(\Lambda_{n}\otimes \IHH)^{n}=\big\{f=(f_{1},\dots,f_{n})^{T}
		\colon f_{j}\in\Lambda_{n}\otimes \IHH \>\>\text{for $j=1,\dots,n$}\big\}
	\end{equation*}
	by setting
	\begin{equation*}
		P^{(n)}\coloneqq
		\begin{pmatrix}
                & & & \hat{\theta}_{n} \otimes P_{n} \\
			\hat{\theta}_{n-1}\otimes P_{n-1} & & & \\
                & \ddots & & \\
			  & & \hat{\theta}_{1}\otimes P_{1} &
		\end{pmatrix},
	\end{equation*}
	on its natural domain of definition given by all $f\in \IHH^{(n)}$ such that $f_{n}\in \Lambda_{n}\otimes \dom(P_{n})$
	and $f_{j}\in \Lambda_{n}\otimes \dom(P_{n-j})$ for all $j=1,\dots,n-1$. Here, for every $\alpha \in \Lambda_{n}$ we use
	$\hat{\alpha}$ to denote the induced endomorphism
	\begin{equation*}
		\hat{\alpha}\colon \Lambda_{n}\longrightarrow\Lambda_{n},\quad \beta\longmapsto \alpha \beta.
	\end{equation*}
	Likewise, we define the selfadjoint operator $H^{(n)}\geq 0$ in the Hilbert space $\IHH^{(n)}$ by the diagonal matrix
	of operators
	\begin{equation*}
		H^{(n)}\coloneqq
		\begin{pmatrix}
			\mathrm{id}_{\Lambda_n}\otimes H & & \\
			& \ddots & \\
			& & \mathrm{id}_{\Lambda_n}\otimes H
		\end{pmatrix}
	\end{equation*}
	on its natural domain of definition given by all $f\in \IHH^{(n)}$ such that $f_{j}\in \Lambda_{n}\otimes \dom(H)$ for
	all $j=1,\dots,n$. With $a\coloneqq\max_{j} a_{j}\in (0,1)$, using
	\begin{equation*}
		P^{(n)}(H^{(n)}+1)^{-a}=
		\begin{pmatrix}
			& & & \hat{\theta}_{n} \otimes P_{n} (H+1)^{-a} \\
			\hat{\theta}_{n-1}\otimes P_{n-1}(H+1)^{-a} & & & \\
			& \ddots & & \\
			& & \hat{\theta}_{1}\otimes P_{1}(H+1)^{-a} &
		\end{pmatrix},
	\end{equation*}
	it follows that $P^{(n)}(H^{(n)}+1)^{-a}$, and hence $P^{(n)}(H^{(n)}+1)^{-1}$, is bounded and
	$\dom(H^{(n)})\subset\dom(P^{(n)})$. The canonical isomorphism
	\begin{equation*}
		\ILL(\IHH^{(n)}) \simeq \Lambda_{n} \otimes \Lambda_{n}^{*} \otimes \ILL(\IHH^{n})
	\end{equation*}
	induces the Fermionic integral
	\begin{equation*}
		\dashint_{ \Lambda_n^* \otimes \ILL(\IHH^n)}\colon \ILL(\IHH^{(n)})\longrightarrow \Lambda_{n}^{*}\otimes \ILL(\IHH^{n}).
	\end{equation*}
	Moreover, we define the linear map
	\begin{equation*}
		\pi^{(n)}_{\IHH}\colon \Lambda_{n}^{*}\otimes \ILL(\IHH^{n})\longrightarrow \ILL( \IHH), \quad \alpha \otimes (a_{ij})_{i,j=1,\dots,n}
		\longmapsto a_{nn}.
	\end{equation*}

	\begin{theorem}\label{main} 
            Under (\ref{ass}), the operator $H^{(n)}+P^{(n)}$ generates an analytic semigroup of angle $\pi/2$ in
		$\IHH^{(n)}$, to be denoted by $\mathrm{e}^{-t(H^{(n)}+P^{(n)})}$, and for all $\geq 0$ one has
		\begin{align}
			\nn & \mathrm{e}^{-t(H^{(n)}+P^{(n)})}=\mathrm{e}^{-tH^{(n)}}+\sum^{n}_{l=1}(-1)^{l}\Phi^{H^{(n)}}_{t}\underbrace{(P^{(n)},\dots,P^{(n)})}_{\text{$l$ times}}, \\
			\label{qayx} & \pi^{(n)}_{\IHH}\dashint_{ \Lambda_n^* \otimes \ILL(\IHH^n)}\mathrm{e}^{-t(H^{(n)}+P^{(n)})}=(-1)^{n}\Phi_{t}^{H}(P_{1},\dots, P_{n}).
		\end{align}
	\end{theorem}

	\begin{proof}
		Step 1: $H^{(n)}+P^{(n)}$ generates an analytic semigroup of angle $\pi/2$.

		Proof of Step 1: As $H^{(n)}$ generates an analytic semigroup of angle $\pi/2$, it suffices to show that $P^{(n)}$
		is $H^{(n)}$-bounded with relative bound $0$ (cf. Theorem 2.10, page 176, in \cite{nagel}). Under (\ref{ass}), this relative
		boundedness can be deduced with some effort using the abstract functional analytic machinery of Favard spaces (cf. Lemma
		2.13, page 176, in \cite{nagel}). We found, nevertheless, an independent and elementary proof that we give here: let
		$f\in \dom(H^{(n)})\subset \dom((H^{(n)}+1)^{a})$, where $a=\max_{j} a_{j}\in (0,1)$. Using the absolute convergence
		of the integral in the formula (cf. Lemma A.4 in \cite{kato3})
		\begin{equation*}
			(H^{(n)}+1)^{a}f=\frac{\sin(\pi a)}{a}\int^{\infty}_{0}s^{a-1}(H^{(n)}+1)( H^{(n)}+1+s)^{-1}f\, \Id s,
		\end{equation*}
		we may estimate for all $\lambda >0$, with some constant $C=C(P^{(n)},H^{(n)},a)>0$, as follows:
		\begin{align*}
                & \left\| P^{(n)}f \right\| \\
			\leq & \left\| P^{(n)}(H^{(n)}+1)^{-a}\right\|\left\|(H^{(n)}+1)^{a}f \right\| \\
			\leq & \; C\int_{0}^{\lambda} s^{a-1}\left\| (H^{(n)}+1)( H^{(n)}+1+s)^{-1}f\right\|\, \Id s+C\int^{\infty}_{\lambda} s^{a-1}\left\| (H^{(n)}+1)( H^{(n)}+1+s)^{-1}f\right\|\, \Id s \\
			\leq & \; C\int^{\lambda}_{0} s^{a-1}\left\| (H^{(n)}+1)( H^{(n)}+1+s)^{-1}\right\|\, \Id s\left\|f\right\| \\
			& +C\int_{\lambda}^{\infty} s^{a-1}\left\| ( H^{(n)}+1+s)^{-1}\right\|\, \Id s \left\|(H^{(n)}+1)f\right\| \\
			\leq & \; C\int^{\lambda}_{0} s^{a-1}\, \Id s\left\|f\right\|+C\int_{\lambda}^{\infty} s^{a-1}s^{-1}\, \Id s \left\|H^{(n)}f\right\|+C\int_{\lambda}^{\infty} s^{a-1}s^{-1}\, \Id s \left\|f\right\|,
		\end{align*}
		where, setting $A=H^{(n)}+1$, we used that for every self-adjoint operator $A\geq 0$ and all $s> 0$ one has
		\begin{equation*}
			\left\|A(A+s)^{-1}\right\|\leq 1, \quad \left\|(A+s)^{-1}\right\|\leq 1/s.
		\end{equation*}
		Thus, by making $\lambda$ large, we get that for all $\varepsilon>0$ there exists $C_{\varepsilon} >0$ such that
		\begin{equation*}
			\left\| P^{(n)}f \right\|\leq \varepsilon \left\|H^{(n)}f\right\|+C_{\varepsilon} \left\| f\right\|\quad\text{for
			all $f\in \dom(H^{(n)})$},
		\end{equation*}
		which means that $P^{(n)}$ is $H^{(n)}$-bounded with relative bound $0$.
		
		That the angle of holomorphy can be chosen maximal follows from the proof of Theorem 2.10, page 176, in \cite{nagel},
		see also page 181 in \cite{nagel}.

		Step 2: For all $t>0$ the series of bounded operators
		\begin{equation*}
			P(t) \coloneqq \mathrm{e}^{-tH^{(n)}}+\sum^{\infty}_{l=1}(-1)^{l}\Phi^{H^{(n)}}_{t}\underbrace{(P^{(n)},\dots,P^{(n)})}
			_{\text{$l$ times}}
		\end{equation*}
		is finite (in fact, cancels for $l\geq n+1$), and one has
		\begin{align}
			\label{uio}P(t)=\mathrm{e}^{-t(H^{(n)}+P^{(n)})}.
		\end{align}

		Proof of Step 2: To see that the sum is finite, we first note that for $s \in [0, \infty)$ and $n \in \IN$ we have
		\begin{equation*}
			P^{(n)}\mathrm{e}^{-sH^{(n)}}=
			\begin{pmatrix}
                    & & & \hat{\theta}_{n} \otimes P_{n} \mathrm{e}^{-sH} \\
				\hat{\theta}_{n-1}\otimes P_{n-1}\mathrm{e}^{-sH} & & & \\
				  & \ddots & & \\
				  & & \hat{\theta}_{1}\otimes P_{1}\mathrm{e}^{-sH} &
			\end{pmatrix}.
		\end{equation*}
		For $s' \in [0, \infty)$ we see that the matrix $P^{(n)}\mathrm{e}^{-sH^{(n)}}P^{(n)}\mathrm{e}^{-s'H^{(n)}}$ has
		precisely one non-zero entry in each column and row which is of the form
		$\hat{\theta}_{i} \hat{\theta}_{j} \otimes P_{i} \mathrm{e}^{-sH}P_{j} \mathrm{e}^{-s'H}$ for some $i \neq j$. Thus,
		if we consider the $l$-th summand in the above series defining $P(t)$ where $l \geq n + 1$, the integrand is a matrix
		for which every non-zero entry is of the form
		\begin{align}
			\label{P^l}\hat{\theta}_{j_1}\cdots \hat{\theta}_{j_l}\otimes \mathrm{e}^{-s_1H}P_{j_{1}}\mathrm{e}^{-(s_2 - s_1)H}P_{j_2}\cdots \mathrm{e}^{-(s_l - s_{l-1})H}P_{j_l}\mathrm{e}^{-(t - s_l)H}
		\end{align}
		for some $1 \leq j_{1}, \ldots, j_{l} \leq n$. Since $l \geq n + 1$, we necessarily have $j_{q} = j_{r}$ for some $1
		\leq q,r \leq l$ so that $\hat{\theta}_{j_1}\cdots \hat{\theta}_{j_l}= 0$. Hence, the sum defining $P(t)$ is finite.

		To see $(\ref{uio})$, by uniqueness of solutions of initial value problems induced by generators of strongly
		continuous semigroups, it suffices to show that for all $f$ in the dense subspace $C^{\infty}(H^{(n)})$ of $\IHH^{(n)}$,
		setting $f(t)\coloneqq P(t)f$, $t>0$, one has
		\begin{itemize}
			\item[(i)] $f(t)\in\dom(H^{(n)})\subset \dom(P^{(n)})$ for all $t>0$,

			\item[(ii)] $f(\bullet)$ is $C^{1}$ on $(0,\infty)$ with
				\begin{equation*}
					(\Id/\Id t) f(t)=-(H^{(n)}+P^{(n)})f(t)\quad\text{for all $t>0$,}
				\end{equation*}

			\item[(iii)] $f(t)\to f$ as $t\to 0+$.
		\end{itemize}
		To this end, for $t\geq 0$, define $f(0;t)\coloneqq \mathrm{e}^{-t H^{(n)}}f$. Clearly, $f(0;\bullet)$ takes values
		in $\dom(H^{(n)})\subset \dom(P^{(n)})$ and
		\begin{equation*}
			P^{(n)}f(0;t)=P^{(n)}(H^{(n)}+1)^{-1}\mathrm{e}^{-t H^{(n)}}(H^{(n)}+1)f
		\end{equation*}
		is $C^{1}$ in $t$, being of the form $A\mathrm{e}^{-t H^{(n)}}g$ with $A$ bounded and $g\in \dom(H^{(n)})$. Assume
		now that for $l\in\IN$ the function
		\begin{equation*}
			f(l;\bullet) \colon [0,\infty)\longrightarrow \IHH^{(n)}
		\end{equation*}
		takes values in $\dom(H^{(n)})\subset \dom(P^{(n)})$ and $P^{(n)}f(l;\bullet)$ is $C^{1}$. Then Lemma \ref{aux} implies
		that
		\begin{equation*}
			t\longmapsto f(l+1;t)\coloneqq-\int^{t}_{0} \mathrm{e}^{-sH^{(n)}}P^{(n)}f(l;t-s)\, \Id s
		\end{equation*}
		is $C^{1}$ with values in $\dom(H^{(n)})\subset \dom(P^{(n)})$ and
		\begin{align}
			\label{ssum}(\Id/\Id t)f(l+1;t)= H^{(n)}f(l+1;t)-P^{(n)}f(l;t).
		\end{align}
		It follows that
		\begin{equation*}
			H^{(n)}f(l+1;t)=P^{(n)}f(l;t)+(\Id/\Id t)f(l+1;t)
		\end{equation*}
		and
		\begin{equation*}
			P^{(n)}f(l+1;\bullet)=P^{(n)}(H^{(n)}+1)^{-1}(H^{(n)}+1)f(l+1;\bullet).
		\end{equation*}
		are $C^{1}$ in $t$. Ultimately, the sequence $f(l;\bullet)$ is well-defined, and Theorem 1.19 on page 488 in
		\cite{kato} also gives
		\begin{align}
			\label{ssum2}\text{$f(l;0)=0$ for all $l\geq 1$, while $f(0;0)=f$. }
		\end{align}
		Finally, by construction we have
		\begin{equation*}
			f(t)=\sum^{\infty}_{l=0}f(l;t)\quad \text{for all $t>0$},
		\end{equation*}
		so (i) is satisfied. The assumption (ii) is satisfied by summing (\ref{ssum}), and (iii) follows from (\ref{ssum2}) and
		the continuity of the $f(l;\cdot)$. This completes the proof of Step 2.

		Step 3: one has (\ref{qayx}).

		Proof of Step 3: By definition of the Fermionic integral and using (\ref{P^l}), we only have to consider the $l = n$
		summand in the series defining $P(t)$. Observe that
		\begin{equation*}
			\begin{pmatrix}
				  & & & 1 \\
				1 & & & \\
				  & \ddots & & \\
				  & & 1 &
			\end{pmatrix}
			\in \mathrm{Mat}_{n\times n}(\IC)
		\end{equation*}
		is the matrix associated to the cyclic permutation
		\begin{equation*}
			\tau \coloneqq
			\begin{pmatrix}
				1 & 2 & \dots & n - 1 & n \\
				2 & 3 & \dots & n & 1
			\end{pmatrix}
		\end{equation*}
		of length $n$. From this it follows that for fixed $(s_{1},\dots, s_{n})\in t\sigma_{n}$ and $1 \leq k \leq n$ the
		matrix
		\begin{equation*}
			\mathrm{e}^{-s_1H}P^{(n)}\mathrm{e}^{-(s_2 - s_1)H}P^{(n)}\dots \mathrm{e}^{-(s_k - s_{k-1})H}P^{(n)}\mathrm{e}^{-(t
			- s_k)H}
		\end{equation*}
		has precisely one non-zero entry in each column and row; if $1 \leq q,r \leq n$, then the $(q,r)$-th entry
		$p^{(n)}(k; q, r)$ is given by
		\begin{equation*}
			\begin{cases}
                    \hat{\theta}_{\tau^{-k + 1}(n-r)}\hat{\theta}_{\tau^{-k + 2}(n-r)}\cdots \hat{\theta}_{n-r}\otimes \mathrm{e}^{-s_1H}P_{\tau^{-k + 1}(n-r)}\mathrm{e}^{-(s_2 - s_1)H}P_{\tau^{-k + 2}(n-r)}\cdots \\
				\hfill \cdots \mathrm{e}^{-(s_k - s_{k-1})H}P_{n-r}\mathrm{e}^{-(t-s_k)H}, & \text{if $r \leq n - 1$,} \\
				& \text{$q = \tau^{k-1}(r + 1)$}; \\
				\hat{\theta}_{\tau^{-k + 1}(n)}\hat{\theta}_{\tau^{-k + 2}(n)}\cdots \hat{\theta}_{n}\otimes \mathrm{e}^{-s_1H}P_{\tau^{-k + 1}(n)}\mathrm{e}^{-(s_2-s_1)H}P_{\tau^{-k + 2}(n)}\cdots \\ 
                    \hfill \cdots \mathrm{e}^{-(s_k - s_{k-1})H}P_{n}\mathrm{e}^{-(t-s_k)H}, & \text{if $r = n$,} \\
                    & \text{$q = \tau^{k-1}(1)$}; \\
				0, & \text{else}.
			\end{cases}
		\end{equation*}
		Using the above formula for $k = n$, it follows that
		\begin{equation*}
			p^{(n)}(n; n, n) = \hat{\theta}_{1} \cdots \hat{\theta}_{n} \otimes \mathrm{e}^{-s_1 H}P_{1} \mathrm{e}^{-(s_2 -
			s_1)H}\cdots \mathrm{e}^{-(s_n - s_{n-1})H}P_{n} \mathrm{e}^{-(t - s_n)H}
		\end{equation*}
		which in turn implies
		\begin{align*}
               \pi^{(n)}_{\IHH}\int_{ \Lambda_n^* \otimes \ILL(\IHH^n)}\mathrm{e}^{-t(H^{(n)}+P^{(n)})} & = (-1)^{n}\int_{t\sigma_n}\mathrm{e}^{-s_1H}P_{1}\mathrm{e}^{-(s_2-s_1)H}P_{2} \times \cdots \\
                & ~~~~~~~~~~~~~~~~~~~~~~~~~ \times \mathrm{e}^{-(s_n-s_{n-1})H}P_{n} \mathrm{e}^{-(t-s_n)H}\, \Id s_{1} \ldots \Id s_{n},
		\end{align*}
		as desired.
	\end{proof}

	\begin{corollary}
		\label{coro} Assume (\ref{ass}).
		
		a) The map
		\begin{equation*}
			[0,\infty)\ni t\longmapsto \Phi^{H}_{t}(P_{1},\dots,P_{n})\in \ILL(\IHH)
		\end{equation*}
		is smooth on $(0,\infty)$ and norm continuous at $t=0$.

		b) For all $t\geq 0$ one has
		\begin{align}\label{smoothing}
                \mathrm{Ran}\big(\Phi^{H}_{t}(P_{1},\dots,P_{n})\big)\subset \dom(H),
		\end{align}
		and for all $t>0$,
		\begin{equation*}
			(\Id/\Id t)\Phi^{H}_{t}(P_{1},\dots,P_{n})= - H \Phi^{H}_{t}(P_{1}, \dots, P_{n}) + P_{1} \Phi^{H}_{t}(P_{2}, \dots
			, P_{n}),
		\end{equation*}
		with the convention that the second summand is $P_{1} \mathrm{e}^{-tH}$ for $n=1$ (in accordance with (\ref{con})).
	\end{corollary}

	\begin{proof}
		We start by remarking that the linear map
		\begin{equation*}
			A\coloneqq\pi^{(n)}_{\IHH}\dashint_{ \Lambda_n^* \otimes \ILL(\IHH^n)}\colon \ILL(\IHH^{(n)})\longrightarrow \ILL(\IHH)
		\end{equation*}
		is continuous, when either both $\ILL(\IHH^{(n)})$ and $\ILL(\IHH)$ are considered as Banach spaces together with
		the operator norm, or both are considered as locally convex spaces equipped with the strong topology.

		a) Clearly, $t\mapsto \mathrm{e}^{-t (H^{(n)}+P^{(n)})}$ is smooth on $(0,\infty)$ (this map even extends
		analytically to the open right complex plane), and
		\begin{equation*}
			\Phi^{H}_{t}(P_{1},\dots,P_{n})=(-1)^{n} A \mathrm{e}^{-t (H^{(n)}+P^{(n)})},
		\end{equation*}
		where $A$ is a continuous and linear, and thus, smooth map between Banach spaces. This shows that $t\mapsto \Phi^{H}_{t}
		(P_{1},\dots,P_{n})$ is smooth, too.
		
		The norm continuity at $t=0$ follows from Lemma \ref{con2}.

		b) For the asserted smoothing property (\ref{smoothing}), note first that by analyticity,
		\begin{align}
			\label{bhu}\mathrm{Ran}\big(\mathrm{e}^{-t (H^{(n)}+P^{(n)})}\big)\subset C^{\infty}(H^{(n)}+P^{(n)})\subset\dom(H^{(n)}+P^{(n)})=\dom(H^{(n)}).
		\end{align}
		Let $f\in \IHH$ and define $\tilde{f}\in \IHH^{(n)}$ by
		\begin{equation*}
			\tilde{f}\coloneqq(0,\dots,0, 1_{\Lambda_n}\otimes f)^{T}\in \IHH^{(n)}=(\Lambda_{n}\otimes \IHH)^{n}.
		\end{equation*}
		Then
		\begin{equation*}
			\Phi^{H}_{t} (P_{1},\dots, P_{n}) f= \quad\text{$n$-th component of
			$(-1)^{n}\dashint_{\IHH^n }\mathrm{e}^{-t (H^{(n)}+P^{(n)})}\tilde{f}$}.
		\end{equation*}
		If we write
		\begin{equation*}
			\mathrm{e}^{-t (H^{(n)}+P^{(n)})}\tilde{f}=\sum^{n}_{k=0}\sum_{1\leq j_1<\ldots <j_k\leq n}\theta_{j_1}\cdots \theta
			_{j_k}\otimes \alpha_{j_1,\ldots,j_k},\quad \alpha_{j_1,\ldots,j_k}\in \IHH^{n},
		\end{equation*}
		then (\ref{bhu}) shows $\alpha_{j_1,\ldots,j_k}\in \dom(H^{(n)})$ for all $j_{1},\ldots,j_{k}$. In particular,
		\begin{equation*}
			\dashint_{\IHH^n }\mathrm{e}^{-t (H^{(n)}+P^{(n)})}\tilde{f}=\alpha_{1,\dots,n}\in \dom(H^{(n)}),
		\end{equation*}
		which proves (\ref{smoothing}).

		For the derivative, since the map $A$ is a bounded operator between Banach spaces, we may calculate
		\begin{align*}
			(\Id/\Id t) \Phi^{H}_{t}(P_{1}, \dots, P_{n}) & = (-1)^{n} (\Id/\Id t) \pi^{(n)}_{\IHH}\dashint_{ \Lambda_n^* \otimes \ILL(\IHH^n)}\mathrm{e}^{-t(H^{(n)} + P^{(n)})} \\
			  & = (-1)^{n} \pi^{(n)}_{\IHH}\dashint_{ \Lambda_n^* \otimes \ILL(\IHH^n)}(\Id/\Id t) \mathrm{e}^{-t(H^{(n)} + P^{(n)})} \\
			& = (-1)^{n+1}\pi^{(n)}_{\IHH}\dashint_{ \Lambda_n^* \otimes \ILL(\IHH^n)}(H^{(n)}+ P^{(n)})\mathrm{e}^{-t(H^{(n)} + P^{(n)})} \\
			  & = -H \Phi^{H}_{t}(P_{1}, \dots, P_{n}) + (-1)^{n+1}\pi^{(n)}_{\IHH}\dashint_{ \Lambda_n^* \otimes \ILL(\IHH^n)}P^{(n)}\mathrm{e}^{-t(H^{(n)} + P^{(n)})}.
		\end{align*}
		Using the series expansion for $\mathrm{e}^{-t(H^{(n)} + P^{(n)})}$ and arguments similar to Step 3 from the proof
		of Theorem 2.1, it follows that
		\begin{equation*}
			\pi^{(n)}_{\IHH}\dashint_{ \Lambda_n^* \otimes \ILL(\IHH^n)}P^{(n)}\mathrm{e}^{-t(H^{(n)} + P^{(n)})}= (-1)^{n-1}P_{1}
			\Phi^{H}_{t}(P_{2}, \dots, P_{n}),
		\end{equation*}
		which implies the claimed formula.
	\end{proof}

	\begin{remark}
		At first sight it seems that the first inclusion of (\ref{bhu}) and the arguments below that inclusion are promising
		enough to conclude that $\Phi^{H}_{t}(P_{1},\dots,P_{n})$ generally maps into $C^{\infty}(H)$. Taking a closer look,
		however, one finds that this is not the case, even for $n=1$: if $a+\theta \otimes b$ is in
		$\dom (H+\hat{\theta}\otimes P)^{2}$, then $a+\theta\otimes b\in \dom(H+\hat{\theta}\otimes P)$, which is equivalent
		to $a,b\in \dom(H+P)=\dom(H)$. Also
		\begin{equation*}
			\dom(H+\hat{\theta}\otimes P)\ni (H+\hat{\theta}\otimes P)(a+\theta \otimes b)= Ha+\theta\otimes (Pa+Hb),
		\end{equation*}
		which is equivalent to $Ha,Pa+Hb\in \dom(H+P)=\dom(H)$. But there is no reason to expect that
		\begin{equation*}
			b=\dashint_{\IHH}(a+\theta\otimes b)
		\end{equation*}
		satisfies $Hb\in \dom(H)$, which would be needed to conclude $b\in \dom(H^{2})$. In fact, Corollary \ref{coro} is optimal in the sense that in general, somewhat surprisingly $\Phi^{H}_{t}(P_{1},\dots
		,P_{n})$ does not even map $C^{\infty}(H)$ into the domain of definition of higher powers of $H$ (even if the
		$P_{j}$'s are bounded). To see this, assume $n=1$, let $P$ be bounded and let $f\in C^{\infty}(H)$. Then for all $t>0$,
		setting $f(r)\coloneqq \mathrm{e}^{-rH}f$, we have
		\begin{equation*}
			\Phi^{H}_{t}(P)f=\int^{t}_{0} \mathrm{e}^{-sH}Pf(t-s) \, \Id s,
		\end{equation*}
		and therefore,
		\begin{equation*}
			(\Id/\Id t)\Phi^{H}_{t}(P)f=\mathrm{e}^{-tH}Pf - \int^{t}_{0} \mathrm{e}^{-sH}PHf(t-s) \, \Id s.
		\end{equation*}
		On the other hand,
		\begin{equation*}
			(\Id/\Id t)\Phi^{H}_{t}(P)f=-H\Phi^{H}_{t}(P)f+Pf(t),
		\end{equation*}
		so that
		\begin{equation*}
			H\Phi^{H}_{t}(P)f=-\mathrm{e}^{-tH}Pf + \int^{t}_{0} \mathrm{e}^{-sH}PHf(t-s) \, \Id s+Pf(t)=-\mathrm{e}^{-tH}Pf +
			\Phi^{H}_{t}(PH)f+P\mathrm{e}^{-tH}f.
		\end{equation*}
		The first and second summand are in $\dom(H)$. Thus $\Phi^{H}_{t}(P)f$ is not in $\dom(H^{2})$ if and only if
		$P\mathrm{e}^{-tH}f$ is not in $\dom(H)$. To construct such an example, let $H$ be the Laplace-Beltrami operator on
		a closed and connected Riemannian manifold $X$, and let $P$ be multiplication by the indicator function $1_{U}$,
		where $U\subset X$ is open with $X\setminus \overline{U}$ nonempty. Then for $f=1\in C^{\infty}(H)=C^{\infty}(X)$ and
		every $t>0$ one has $\mathrm{e}^{-tH}1=1$ as $X$ is closed. Consequently, $P\mathrm{e}^{-tH}f=1_{U}$, which is not
		in $\dom(H)=W^{2}(X)$. In fact, $1_{U}$ is not even in $W^{1}(X)$.
	\end{remark}

	\section{Application to Geometric Operators}

	Assume $X$ is a closed connected Riemannian manifold with its volume measure $\mu$ and geodesic distance $d(x,y)$, let
	$E\to X$ be an Hermitian vector bundle, $H$ is a formally self-adjoint, elliptic, nonnegative differential operator of
	order $2$ on $E\to X$, and $P_{1},\dots, P_{n}$ are differential operators on $E\to X$ of order $\leq 1$. We consider $H$
	as a self-adjoint operator in the Hilbert space $\Gamma_{L^2}(X,E)$ of $L^{2}$-sections with domain of definition the Sobolev
	space $\Gamma_{W^2}(X,E)$, and each $P_{j}$ as a closed densely defined operator in $\Gamma_{L^2}(X,E)$ with domain of
	definition $\Gamma_{W^1}(X,E)$. Then $(H+1)^{-1/2}$ is a pseudo-differential operator of order $-1$ which implies that
	$P_{j} (H+1)^{-1/2}$ is a pseudo-differential operator of order $0$ and thus bounded in $\Gamma_{L^2}(X,E)$. Let us
	also record that $H^{(n)}$ is a self-adjoint elliptic differential operator on the Hermitian vector bundle
	\begin{equation*}
		\Lambda_{n}\otimes E\otimes \IC^{n}\longrightarrow X,
	\end{equation*}
	and $P^{(n)}$ is a differential operator on the latter bundle of order $\leq 1$. Let
	\begin{equation*}
		E\boxtimes E^{*}\longrightarrow X\times X
	\end{equation*}
	denote the smooth vector bundle whose fiber $(x,y)$ is given by $\Hom(E_{y},E_{x})$ and let $X\bowtie X$ be the open subset
	of $X\times X$ given by
	\begin{equation*}
		X\bowtie X=\{(x,y)\in X\times X \colon \text{there is a unique minimizing geodesic connecting $x$ and $y$}\}.
	\end{equation*}
	With these preparations, we can now prove:

	\begin{theorem}
		\label{expan}\emph{a)} $\Phi^{H}(P_{1},\dots,P_{n})$ has a smooth integral kernel, in the sense that there exists a
		(uniquely determined) smooth map
		\begin{equation*}
			(0,\infty)\times X\times X \ni (t,x,y)\longmapsto \Phi^{H}_{t}(P_{1},\dots,P_{n})(x,y)\in \Hom(E_{y},E_{x})\subset
			E\boxtimes E^{*},
		\end{equation*}
		such that for all $t>0$, $f\in \Gamma_{L^2}(X,E)$, $\mu$-a.e. $x\in X$,
		\begin{equation*}
			\Phi^{H}_{t}(P_{1},\dots,P_{n})f(x)=\int_{X} \Phi^{H}_{t}(P_{1},\dots,P_{n})(x,y)f(y) \, \Id\mu(y).
		\end{equation*}
		\emph{b)} If $H$ has a scalar principal symbol, then the integral kernel $\Phi^{H}_{t}(P_{1},\dots,P_{n})(x,y)$ has
		a strong asymptotic expansion as $t\to 0$+, meaning that there exists a (uniquely determined) sequence of smooth
		maps
		\begin{equation*}
			(0,\infty)\times X\bowtie X \ni (t,x,y)\longmapsto \Phi^{H}_{t,j}(P_{1},\dots,P_{n})(x,y)\in \Hom(E_{y},E_{x})\subset
			E\boxtimes E^{*},\quad j\in\IN_{\geq 0},
		\end{equation*}
		with the following property: for every Hermitian connection $\nabla$ on $E\to X$, every $K\subset X\bowtie X$
		compact, and every $T>0$, $l, k_{1},k_{2},k_{3}\in\IN_{\geq 0}$, there exists a constant $C>0$, such that for all
		$0<t\leq T$, $(x,y)\in K$,
		\begin{align*}
                \Big|\partial_{t}^{k_1}\nabla_{x}^{k_2}\nabla_{y}^{k_3}\Big((4\pi t)^{\frac{\dim(X)}{2}}\mathrm{e}^{\frac{d(x,y)^{2}}{4t}}\Phi^{H}_{t}(P_{1},\dots,P_{n})(x,y)-\sum^{l}_{j=0}\frac{t^{j}}{j!}\Phi^{H}_{t,j}(P_{1},\dots,P_{n}     (x,y)\Big)& \Big|_{(x,y)}  \\
                & \leq C t^{l+1-k_1}.
		\end{align*}
	\end{theorem}

	\begin{proof}
		a) Note first that by local parabolic regularity, $H^{(n)}+P^{(n)}$ has a smooth heat kernel
		\begin{align*}
			(0,\infty)\times X\times X\ni (t,x,y)\longmapsto \mathrm{e}^{-t (H^{(n)}+P^{(n)})}(x,y)\in & \>\>\Lambda_{n}\otimes \Lambda_{n}^{*}\otimes \mathrm{Mat}_{n\times n}(\Hom(E_{y},E_{x})) \\
			\ & \subset (\Lambda_{n} \otimes E \otimes \IC^{n})\boxtimes (\Lambda_{n} \otimes E \otimes \IC^{n})^{*}.
		\end{align*}
		Let
		\begin{align}
			\label{utu}A \colon \Lambda_{n}\otimes \Lambda_{n}^{*}\otimes \mathrm{Mat}_{n\times n}(E\boxtimes E^{*})) \longrightarrow E\boxtimes E^{*}
		\end{align}
		denote the homomorphism of smooth vector bundles over $X\times X$, given at $(x,y)$ by the composition of the linear
		maps
		\begin{equation*}
			\dashint_{\Lambda^*_n\otimes \mathrm{Mat}_{n\times n}(\Hom(E_y,E_x))}\colon \Lambda_{n}\otimes \Lambda_{n}^{*}
			\otimes \mathrm{Mat}_{n\times n}(\Hom(E_{y},E_{x})) \longrightarrow\Lambda_{n}^{*}\otimes \mathrm{Mat}_{n\times n}(
			\Hom(E_{y},E_{x}))
		\end{equation*}
		and
		\begin{equation*}
			\Lambda_{n}^{*}\otimes \mathrm{Mat}_{n\times n}(\Hom(E_{y},E_{x}))\longrightarrow \Hom(E_{y},E_{x}),\quad \alpha \otimes
			(a_{ij})_{i,j=1,\dots,n}\mapsto (-1)^{n}a_{nn}.
		\end{equation*}
		Then, in view of Theorem \ref{main}, the assignment
		\begin{align}
			\label{utu2}\Phi^{H}_{t}(P_{1},\dots,P_{n})(x,y)\coloneqq A(x,y)\mathrm{e}^{-t (H^{(n)}+P^{(n)})}(x,y)
		\end{align}
		does the job.
		
		b) As $H^{(n)}+P^{(n)}$ has a scalar symbol, it follows from Theorem 1.1 in \cite{lude} that its heat kernel has a
		strong asymptotic expansion, given by, say
		\begin{align*}
			(0,\infty)\times X\times X\ni (t,x,y)\longmapsto \tilde{\Phi}_{t,j}(x,y)\in & \>\>\Lambda_{n}\otimes \Lambda_{n}^{*}\otimes \mathrm{Mat}_{n\times n}(\Hom(E_{y},E_{x})) \\
			\ & \subset (\Lambda_{n} \otimes E \otimes \IC^{n})\boxtimes (\Lambda_{n} \otimes E \otimes \IC^{n})^{*}, \quad j\in \IN_{\geq 0}.
		\end{align*}
		Then, in view of the Leibniz formula for differentiation, we may set
		\begin{equation*}
			\Phi^{H}_{t,j}(P_{1},\dots,P_{n})(x,y)\coloneqq A(x,y)\tilde{\Phi}_{t,j}(x,y),
		\end{equation*}
		completing the proof.
	\end{proof}

	Let us now assume that $H$ is of the form $\nabla^{\dagger}\nabla/2+W$, where $\nabla$ is an Hermitian connection on
	$E\to X$ and $W\in \Gamma_{C^\infty}(X,\End(E))$ is self-adjoint. We decompose every first-order differential
	operator $P$ on $E\to X$ of order $\leq 1$ into its first-order and its zeroth-order part
	\begin{equation*}
		P_{j}= \mathrm{Sym}^{(1)}(P_{j})\circ\nabla +V^{\nabla}_{P_j}
	\end{equation*}
	with respect to $\nabla$, where
	\begin{equation*}
		\mathrm{Sym}^{(1)}(P_{j})\in \Gamma_{C^\infty}(X,\Hom(T^{*}X,\End(E)))
	\end{equation*}
	denotes the first-order symbol of $P_{j}$, and
	\begin{equation*}
		V^{\nabla}_{P_j}\coloneqq P_{j}- \mathrm{Sym}^{(1)}(P_{j})\circ\nabla \in \Gamma_{C^\infty}(X,\End(E)).
	\end{equation*}
	has order $0$.
	
	We fix a filtered probability space with a right-continuous and locally complete filtration which carries for every $x,
	y\in X$, $t>0$ an adapted Brownian bridge $\mathsf{B}^{x,y;t}$ in $X$ from $x$ to $y$ defined on $[0,t]$ (cf. Section \ref{appe}). It follows that
	$\mathsf{B}^{x,y;t}$ is a continuous semimartingale which satisfies
	\begin{equation*}
		\IP\{\mathsf{B}^{x,y;t}_{0}=x, \mathsf{B}^{x,y;t}_{t}=y\} = 1.
	\end{equation*}
        In the sequel, $*\delta$ stands for Stratonovic integration and $\delta$ will denote It\^{o} integration. As for any continuous
	semimartingale on a manifold, one can naturally define the stochastic parallel transport with respect to $\nabla$
	along $\mathsf{B}^{x,y;t}$. It is the uniquely determined continuous semimartingale in $E\boxtimes E^{*}$ defined on $[
	0,t]$ such that \cite{norris} for all $s\in [0,t]$,
	\begin{itemize}
		\item one has $\transport_{\nabla}^{x,y;t}(s)\in \Hom(E_{x},E_{\mathsf{B}_s^{x,y;t}})$ unitarily,

		\item for all $f\in \Gamma_{C^{\infty}}(X,E)$,
			\begin{align}
				\label{para}\transport_{\nabla}^{x,y;t}(s)^{-1}f(\mathsf{B}^{x,y;t}_{s})= \int^{s}_{0}\transport_{\nabla}^{x,y;t}(r)^{-1}\nabla f(*\delta \mathsf{B}_{r}^{x,y;t}),\quad \transport_{\nabla}^{x,y;t}(0)=\mathrm{id}_{E_x},
			\end{align}
	\end{itemize}
	noting that the Stratonovic line integral on the RHS can be defined in coordinates \cite{norris}. Let the process
	\begin{equation*}
		\mathscr{W}_{\nabla}^{x,y;t}(s)\in \End(E_{x}),\quad\text{$s\in [0,t]$,}
	\end{equation*}
	be given as the solution of the pathwise ordinary differential equation
	\begin{equation*}
		(\Id/\Id s) \mathscr{W}_{\nabla}^{x,y;t}(s)=- \mathscr{W}_{\nabla}^{x,y;t}(s)\pa_{\nabla}^{x,y;t}(s)^{-1}W(\mathsf{B}
		^{x,y;t}_{s})\pa_{\nabla}^{x,y;t}(s),\quad \mathscr{W}_{\nabla}^{x,y;t}(0)=\mathrm{id}_{E_x}.
	\end{equation*}

	The pure first-order parts give rise to the It\^{o} line integrals
	\begin{align*}
		\int^{s}_{0}\pa_{\nabla}^{x,y;t}(r)^{-1}\mathrm{Sym}^{(1)}(P_{j})^{\flat}(\delta \mathsf{B}^{x,y;t}_{r})\pa_{\nabla}^{x,y;t}(r)\in \End(E_{x}),\quad\text{$s\in [0,t]$},
	\end{align*}
	noting that the Itô differential $\delta \mathsf{B}^{x,y;t}$ is understood with respect to the Levi-Civita connection on $TX$. From the zeroth-order parts, on the other hand, we get upon pathwise Riemann integration the processes
	\begin{equation*}
		\int^{s}_{0}\pa_{\nabla}^{x,y;t}(r)^{-1}V_{P_j}^{\nabla}(\mathsf{B}^{x,y;t}_{r})\pa_{\nabla}^{x}(r) \, \Id r \in \End
		(E_{x}),\quad\text{$s\in [0,t]$,}
	\end{equation*}
	which are adapted locally absolutely continuous processes. We can then define a continuous semimartingale
	\begin{equation*}
		\int_{t\sigma_n}\delta\Psi^{x,y;t}_{\nabla,P_1}(s_{1})\cdots \delta\Psi^{x,y;t}_{\nabla,P_n}(s_{n}) \colon [0,\infty)
		\times \Omega\longrightarrow \End(E_{x}),
	\end{equation*}
	as follows: with
	\begin{equation*}
		\Psi^{x,y;t}_{\nabla,P_i}(s) \coloneqq \int^{s}_{0}\pa_{\nabla}^{x,y;t}(r)^{-1}\big(\mathrm{Sym}^{(1)}(P_{i})^{\flat}
		(\delta\mathsf{B}^{x,y;t}_{r}) + V^{\nabla}_{P_i}(\mathsf{B}^{x,y;t}_{r})\big)\pa_{\nabla}^{x}(r) \, \Id r \in \End
		(E_{x}),\quad\text{$s\in [0,t]$},
	\end{equation*}
	we define an iterated It\^{o} integral inductively by
	\begin{align*}
            & \int_{s\sigma_1}\delta\Psi^{x,y;t}_{\nabla,P_1}(s_{1})\coloneqq\Psi^{x,y;t}_{\nabla,P_1}(s), \\
            & \int_{s\sigma_{n+1}}\delta\Psi^{x,y;t}_{\nabla,P_1}(s_{1})\cdots\delta\Psi^{x,y;t}_{\nabla,P_{n+1}}(s_{n+1})\coloneqq\int^{s}_{0}\left(\int_{r\sigma_n}\delta\Psi^{x,y;t}_{\nabla,P_1}(r_{1})\cdots \delta\Psi^{x,y;t}_{\nabla,P_n}(r_{n})\right)\delta\Psi^{x,y;t}_{\nabla,P_{n+1}}(r).
	\end{align*}
	Let
	\begin{equation*}
		(0,\infty)\times X \times X\ni (t,x,y)\longmapsto p(t,x,y)\coloneqq\mathrm{e}^{-\frac{t}{2}\Delta}(x,y)\in (0,\infty)
	\end{equation*}
	be the heat kernel of $-\Delta/2$, noting that $-\Delta/2\geq 0$ with domain of definition $W^{2,2}(X)$ is self-adjoint
	in $L^{2}(X)$, in view of the compactness of $X$. Here comes the main result of this section:

	\begin{theorem}
		\label{Feynman-Kac} For all $t>0$, $x,y\in X$,
		\begin{align*}
			\Phi^{\nabla^\dagger\nabla/2 +W}_{t}(P_{1},\dots,P_{n})(x,y)= p(t,x,y)\mathbb{E}\Big[ \mathscr{W}_{\nabla}^{x,y;t}(t) \int_{t\sigma_n}\delta\Psi^{x,y;t}_{\nabla,P_1}(s_{1})\cdots \delta\Psi^{x,y;t}_{\nabla,P_n}(s_{n}) \pa^{x,y;t}_{\nabla}(t)^{-1}\Big].
		\end{align*}
	\end{theorem}

	\begin{proof}
		Note first that with $H\coloneqq \nabla^{\dagger}\nabla/2 +W$, the operator $H^{(n)}+P^{(n)}$ is given by
		\begin{equation*}
			H^{(n)}+P^{(n)}=\tilde{\nabla}^{\dagger}\tilde{\nabla} / 2 + \tilde{Q},
		\end{equation*}
		where $\tilde{\nabla}$ is the Hermitian connection on
		\begin{equation*}
			\tilde{E}\coloneqq \Lambda_{n}\otimes E\otimes \IC^{n}\longrightarrow X
		\end{equation*}
		given by
		\begin{equation*}
			\tilde{\nabla}\coloneqq
			\begin{pmatrix}
				\mathrm{id}_{\Lambda_n}\otimes\nabla & & \\
                    & \ddots & \\
				& & \mathrm{id}_{\Lambda_n}\otimes \nabla
			\end{pmatrix}
		\end{equation*}

		and $\tilde{Q}$ is an operator of order $\leq 1$ on $\tilde{E}\to X$, given by
		\begin{equation*}
			\tilde{Q}=
			\begin{pmatrix}
				\mathrm{id}_{\Lambda_n}\otimes W & &  \\
				& \ddots & \\
				& & \mathrm{id}_{\Lambda_n}\otimes W
			\end{pmatrix}+
			\begin{pmatrix}
                    & &  & \hat{\theta}_{n} \otimes P_{n} \\
				\hat{\theta}_{n-1}\otimes P_{n-1} &  &  & \\
				  & \ddots & & \\
				  & & \hat{\theta}_{1}\otimes P_{1} &
			\end{pmatrix}.
		\end{equation*}
		It follows that
		\begin{equation*}
			\tilde{Q}=\mathrm{Sym}^{(1)}(\tilde{Q})\tilde{\nabla}+\tilde{W}+\tilde{V},
		\end{equation*}
		where the first-order symbol $\mathrm{Sym}^{(1)}(\tilde{Q})$ is given by
		\begin{equation*}
			\mathrm{Sym}^{(1)}(\tilde{Q})=
			\begin{pmatrix}
                    & & & \hat{\theta}_{n} \otimes \mathrm{Sym}^{(1)}(P_{n}) \\
				\hat{\theta}_{n-1}\otimes \mathrm{Sym}^{(1)}(P_{n-1}) & & & \\
				& \ddots & & \\
				  & & \hat{\theta}_{1}\otimes \mathrm{Sym}^{(1)}(P_{1}) &
			\end{pmatrix}
		\end{equation*}
		and the zeroth-order part is given by
		\begin{align*}
                & \tilde{W}+\tilde{V}\coloneqq 
                \begin{pmatrix}
                    \mathrm{id}_{\Lambda_n}\otimes W & & \\
                    & \ddots &\\
                    & & \mathrm{id}_{\Lambda_n}\otimes W
                \end{pmatrix} + 
                \begin{pmatrix}
                    & & & \hat{\theta}_{n} \otimes V^{\nabla}_{P_n} \\ 
                    \hat{\theta}_{n-1}\otimes V^{\nabla}_{P_{n-1}} & & & \\
                    & \ddots & & \\
                    & & \hat{\theta}_{1}\otimes V^{\nabla}_{P_1} &
                \end{pmatrix}.
		\end{align*}

		According to the main result of \cite{boldt} one has
		\begin{align}
			\label{fey}\mathrm{e}^{-t (\tilde{\nabla}^\dagger \tilde{\nabla}/2+\tilde{Q})}(x,y)= p(t,x,y) \mathbb{E}\left[\mathcal{\tilde{Q}}^{x,y;t}(t)\pa_{\tilde{\nabla}}^{x,y;t}(t)^{-1}\right],
		\end{align}
		where
		\begin{equation*}
			\mathcal{\tilde{Q}}^{x,y;t}(s) \in \End(\tilde{E}_{x}),\quad\text{$s\in [0,t]$,}
		\end{equation*}
		is the solution of the It\^{o} equation
		\begin{equation*}
			\delta \mathcal{\tilde{Q}}^{x,y;t}(s)=- \mathcal{\tilde{Q}}^{x,y;t}(s)\big(\delta\tilde{A}^{x,y;t}_{s}+\delta\tilde
			{B}^{x,y;t}_{s}\big),\quad \mathcal{\tilde{Q}}^{x,y;t}(0)=\mathrm{id}_{\tilde{E}_x},
		\end{equation*}
		with
		\begin{align*}
                & \tilde{A}^{x,y;t}_{s}\coloneqq \int^{s}_{0}\pa_{\tilde{\nabla}}^{x}(r)^{-1}\mathrm{Sym}^{(1)}(\tilde{Q})^{\flat}(\delta \mathsf{B}^{x,y;t}_{r})\pa_{\tilde{\nabla}}^{x,y;t}(r) \in \End(\tilde{E}_{x}), \\
                & \tilde{B}^{x,y;t}_{s}\coloneqq \int^{s}_{0}\pa_{\tilde{\nabla}}^{x,y;t}(r)^{-1}(\tilde{W}+\tilde{V})(\mathsf{B}^{x,y;t}_{r})\pa_{\tilde{\nabla}}^{x,y;t}(r) \, \Id r \in \End(\tilde{E}_{x}),\quad\text{$s\in [0,t]$.}
		\end{align*}
		In fact, by decomposing the zeroth-order contributions, the It\^{o} product rule gives
		\begin{equation}
			\label{Splitting P}\mathcal{\tilde{Q}}^{x,y;t}(s)=\mathscr{\tilde{W}}^{x,y;t}(s)\mathcal{\tilde{Q}}^{x,y;t}_{\#}(s)
			,
		\end{equation}
		where the processes
		\begin{equation*}
			\mathscr{\tilde{W}}^{x,y;t}(s),\mathcal{\tilde{Q}}^{x,y;t}_{\#}(s)\in \End(\tilde{E}_{x}),\quad\text{$s\in
			[0,t]$,}
		\end{equation*}
		are respectively defined by
		\begin{equation*}
			(\Id/\Id s) \mathscr{\tilde{W}}^{x,y;t}(s)=- \mathscr{\tilde{W}}^{x,y;t}(s)\pa_{\tilde{\nabla}}^{x,y;t}(s)^{-1}\tilde
			{W}(\mathsf{B}^{x,y;t}_{s})\pa_{\tilde{\nabla}}^{x,y;t}(s),\quad \mathscr{\tilde{W}}^{x,y;t}(0)=\mathrm{id}_{\tilde{E}_x}
			,
		\end{equation*}
		and
		\begin{equation*}
			\delta \mathcal{\tilde{Q}}^{x,y;t}_{\#}(s)= - \mathcal{\tilde{Q}}^{x,y;t}_{\#}(s) \big(\delta\tilde{A}^{x,y;t}_{s}
			+ \delta\tilde{C}^{x,y;t}_{s}\big),\quad \mathcal{\tilde{Q}}^{x,y;t}_{\#}(0)=\mathrm{id}_{\tilde{E}_x},
		\end{equation*}
		with
		\begin{equation*}
			\tilde{C}^{x,y;t}_{s}\coloneqq \int^{s}_{0}\pa_{\tilde{\nabla}}^{x,y;t}(r)^{-1}\tilde{V}(\mathsf{B}^{x,y;t}_{r})\pa
			_{\tilde{\nabla}}^{x,y;t}(r) \, \Id r \in \End(\tilde{E}_{x}).
		\end{equation*}
		Note that a priori one has the convergence in probability
		\begin{align*}
                & \mathcal{\tilde{Q}}^{x,y;t}_{\#}(t)=\mathrm{id}_{\tilde{E}_x}+\sum^{\infty}_{l=1}(-1)^{l}\int_{t\sigma_l}\big(\delta\tilde{A}^{x,y;t}_{s_1}+ \delta\tilde{C}^{x,y;t}_{s_1}\big)\cdots \big(\delta\tilde{A}^{x,y;t}_{s_l}+ \delta\tilde{C}^{x,y;t}_{s_l}\big),
		\end{align*}
		as It\^{o} equations can be solved via Picard iteration, while a posteriori the sum cancels for $l>n$. Recalling the
		definition of $\Psi^{x,y;t}_{\nabla, P_i}$ and using
		\begin{align}
			\label{epao}\pa_{\tilde{\nabla}}^{x,y;t}(s)= \begin{pmatrix}\mathrm{id}_{\Lambda_n}\otimes\pa_{\nabla}^{x,y;t}(s)&&\\&\ddots&\\&&\mathrm{id}_{\Lambda_n}\otimes \pa_{\nabla}^{x,y;t}(s)\end{pmatrix},
		\end{align}
		we have
		\begin{equation*}
			\delta\tilde{A}^{x,y;t}_{s} + \delta\tilde{C}^{x,y;t}_{s} =
			\begin{pmatrix}
                    & & & \hat{\theta}_{n} \otimes \delta\Psi^{x,y;t}_{\nabla,P_n}(s) \\
				\hat{\theta}_{n-1}\otimes \delta\Psi^{x,y;t}_{\nabla,P_{n-1}}(s) & & & \\
				& \ddots & & \\
				& & \hat{\theta}_{1}\otimes \delta\Psi^{x,y;t}_{\nabla,P_1} &
			\end{pmatrix}.
		\end{equation*}
		The calculations from the proof of Step 3 in Theorem \ref{main} and \eqref{Splitting P} yield
		\begin{align*}
                & \int_{t \sigma_n}\big(\delta\tilde{A}^{x,y;t}_{s_1}+ \delta\tilde{C}^{x,y;t}_{s_1}\big) \cdots\big(\delta\tilde{A}^{x,y;t}_{s_n}+\delta\tilde{C}^{x,y;t}_{s_n}\big) \\
                & ~~~~~~~~~~ = 
                \begin{pmatrix}
                    \ast & & & \\
                    & \ddots & &\\
                    & & \ast &\\
                    & & & \hat{\theta}_{1} \cdots \hat{\theta}_{n} \otimes \int_{t \sigma_n}\delta\Psi_{\nabla,P_1}^{x,y;t}(s_{1}) \dots \delta\Psi_{\nabla,P_n}^{x,y;t}(s_{n})
                \end{pmatrix} .
		\end{align*}
		With
		\begin{equation*}
			\mathscr{\tilde{W}}^{x,y;t}(s)=
			\begin{pmatrix}
				\mathrm{id}_{\Lambda_n}\otimes \mathscr{W}_{\nabla}^{x,y;t}(s) & & \\
                    & \ddots & \\
                    & & \mathrm{id}_{\Lambda_n}\otimes \mathscr{W}_{\nabla}^{x,y;t}(s)
			\end{pmatrix}
		\end{equation*}
		and with $A(x,y)$ as in the proof of Theorem \ref{expan}, using (\ref{epao}) we thus get
		\begin{align*}
                & A(x,y) \left(\mathcal{\tilde{Q}}^{x,y;t}(t)\pa_{\tilde{\nabla}}^{x,y;t}(t)^{-1}\right) \\
                & = \mathscr{W}_{\nabla}^{x,y;t}(t) \, A(x,y) \left(\mathcal{\tilde{Q}}^{x,y;t}_{\#}(t) \pa_{\tilde{\nabla}}^{x,y;t}(t)^{-1}\right) \\
                & = (-1)^{n} \, \mathscr{W}_{\nabla}^{x,y;t}(t) A(x,y) \left( \int_{t \sigma_n}\big(\delta\tilde{A}^{x,y;t}_{s_1}+ \delta\tilde{C}^{x,y;t}_{s_1}\big) \cdots\big(\delta \tilde{A}^{x,y;t}_{s_n}+\delta\tilde{C}^{x,y;t}_{s_n}\big)\pa_{\tilde{\nabla}}^{x,y;t}(t)^{-1}\right) \\
                & = (-1)^{2n}\, \mathscr{W}_{\nabla}^{x,y;t}(t) \int_{t \sigma_n}\delta\Psi_{\nabla,P_1}^{x,y;t}(s_{1}) \dots \delta\Psi_{\nabla,P_n}^{x,y;t}(s_{n})\pa_{\nabla}^{x,y;t}(t)^{-1},
		\end{align*}
		so that the proof is completed using
		\begin{align*}
			\Phi^{\nabla^\dagger\nabla/2 +W}_{t}(P_{1},\dots,P_{n})(x,y)=A(x,y)\mathrm{e}^{ -t (\tilde{\nabla}^\dagger \tilde{\nabla}/2+\tilde{Q})}(x,y)
		\end{align*}
		and (\ref{fey}).
	\end{proof}

	\section{Stochastic Localization on Loop Spaces }

	Assume now that $X$ is a closed and even-dimensional Riemannian spin manifold. The Hermitian vector bundle
	$\Sigma\to X$ (the spinor bundle) is equipped with the spin Levi-Civita connection $\nabla$ and Clifford
	multiplication
	\begin{equation*}
		\mathbf{c}\colon \Omega(X) \longrightarrow \Gamma_{C^\infty}(X,\End(\Sigma)).
	\end{equation*}
	The bundle $\Sigma\to X$, and therefore $\End(\Sigma)\to X$, is fiberwise $\IZ_{2}$-graded and we obtain a fiberwise defined supertrace
	\begin{equation*}
		\Str\colon\Gamma_{C^\infty}(X,\End(\Sigma))\longrightarrow C^{\infty}(X).
	\end{equation*}
	Clifford multiplication becomes an even isomorphism, where the exterior bundle is given the even/odd $\IZ_{2}$-grading.

	Let $D$ denote the spin Dirac operator in $\Sigma\to X$. The closure of this operator (to be denoted with the same symbol)
	induces a self-adjoint and odd operator in $\Gamma_{L^2}(X, \Sigma)$ with domain of definition $W^{1,2}(X,\Sigma)$. Thus,
	$D^{2}$ becomes a self-adjoint nonnegative operator in $\Gamma_{L^2}(X, \Sigma)$ with domain of definition
	$W^{2,2}(X,\Sigma)$.

	Let $u$ be a formal variable of degree $-1$ with $u^{2}=0$ and consider the (locally convex and unital)
	differential graded algebra $\Omega_{\mathbb{T}}(X)$, given by all $\omega = \omega' + u \omega''$ with $\omega', \omega
	'' \in \Omega(X)$. Given such $\omega, \omega_{1}, \ldots, \omega_{n} \in \Omega_{\mathbb{T}}(X)$, we may define differential
	operators
	\begin{align}\label{Components P}
		\begin{split}
                P(\omega)&\coloneqq [D, \mathbf{c}(\omega')] - \mathbf{c}(\Id \omega') + \mathbf{c}(\omega'') \\
                P(\omega_{1}, \omega_{2})&\coloneqq (-1)^{\deg(\omega_1')}\big(\mathbf{c}(\omega_{1}' \wedge \omega_{2}') - \mathbf{c}(\omega_{1}') \mathbf{c}(\omega_{2}')\big) \\ P(\omega_{1}, \ldots, \omega_{k})&\coloneqq 0, \quad k \geq 3
		\end{split}
	\end{align}
	on $\Sigma\to X$ of degrees $\leq 1$, $0$ and $0$ respectively. Denote by $\mathcal{P}_{m, n}$ the set of ordered partitions
	of length $m$ of the set $\{1, \ldots, n\}$, i. e. $m$-tuples of nonempty subsets $I=(I_{1}, \ldots, I_{m})$ of $\{1,\ldots
	,n\}$ such that $I_{1} \cup \dots \cup I_{m} = \{1, \dots, n\}$ and such that each element of $I_{a}$ is smaller than
	any element of $I_{b}$ whenever $a < b$. Given $\omega_{0}, \ldots, \omega_{n} \in \Omega_{\mathbb{T}}(X)$ and
	$I\in \mathcal{P}_{m, n}$, we define
	\begin{equation*}
		\text{$\omega_{I_a}\coloneqq(\omega_{i+1},\dots, \omega_{i+i'})$, if $I_{a}= \{j \colon i < j \leq i + i'\}$ for some $i$,
		$i'$},
	\end{equation*}
	so that each $P(\omega_{I_a})$ is defined according to the above recipe.

	Define a based Brownian loop $\mathsf{B}^{x;t}\coloneqq\mathsf{B}^{x,x;t}$ and denote its associated parallel transport
	with respect to the spin Levi-Civita connection with $\pa^{x;t}(\cdot)$. We also set
	\begin{equation*}
		\Psi^{x;t}_{\omega_{I_a}}\coloneqq\Psi^{x,x;t}_{P(\omega_{I_a})}.
	\end{equation*}

	Here comes the main result of this section:

	\begin{theorem}
		\label{local} Given $\omega_{0}, \ldots, \omega_{n} \in \Omega_{\mathbb{T}}(X)$, $x\in M$, $t>0$, define a complex
		valued random variable by
		\begin{align}
                F_{x;t}(\omega_{0},\dots,\omega_{n}) \coloneqq & \>(t/2)^{ -\frac{n}{2}+\frac{1}{2}\sum^n_{j=0}\mathrm{deg}(\omega'_j) }\Str_{x} \Big(\sum_{m = 1}^{n} (-2)^{m} \sum_{I \in \mathcal{P}_{m,n}}\mathbf{c}_{x}(\omega_{0}')p(t,x,x) \> \times \\
			\nn & \times \mathrm{e}^{-\frac{1}{8}\int_0^t \mathrm{scal}(\mathsf{B}_s^{x;t}) \, \Id s}\int_{t\sigma_n}\delta\Psi^{x;t}_{ \omega_{I_1}}(s_{1})\cdots \delta\Psi^{x;t}_{ \omega_{I_m}}(s_{m}) \pa^{x;t}_{\nabla}(t)^{-1}\Big),
		\end{align}
		with the convention
		\begin{equation*}
			F_{x;t}(\omega_{0})\coloneqq(t/2)^{ \frac{\mathrm{deg}(\omega'_{0})}{2} }\Str_{x} \Big(\mathbf{c}_{x}(\omega_{0}')
			p(t,x,x) \mathrm{e}^{-\frac{1}{8}\int_0^t \mathrm{scal}(\mathsf{B}_s^{x;t}) \, \Id s}\pa^{x;t}_{\nabla}(t)^{-1}\Big).
		\end{equation*}
		Then, uniformly in $x\in X$, the limit
		\begin{align}
			\label{Limit Str}F_{x}(\omega_{0},\dots,\omega_{n}) \coloneqq\lim_{t\to 0+}F_{x;t}(\omega_{0},\dots,\omega_{n})\quad\text{exists in $\bigcap_{b\in [1,\infty)}L^{b}(\IP)$,}
		\end{align}
		and the complex valued random field $F_{\bullet}(\omega_{0},\dots,\omega_{n})$ on $X$ satisfies
		\begin{equation*}
			\mathbb{E}\left[F_{\bullet}(\omega_{0},\dots,\omega_{n})\right]\cdot \underline{\mu}= \frac{(-1)^{n} \, 2^{2n}}{n!
			(2 \pi \sqrt{-1})^{\frac{\dim(X)}{2}}}\big(\hat{A}(X)\wedge \omega'_{0}\wedge\omega''_{1}\wedge\cdots \wedge \omega
			_{n}''\big)^{[d]},
		\end{equation*}
		where $\underline{\mu}$ denotes the Riemannian volume form on $X$ and $(\dots)^{[d]}$ the top-degree part of the
		inhomogeneous differential form.
	\end{theorem}

	In view of the Lichnerowicz formula
	\begin{equation*}
		D^{2}/2 = \nabla^{\dagger}\nabla/2 + \mathrm{scal}/8
	\end{equation*}
	and Theorem \ref{Feynman-Kac}, we get the fundamental relation
	\begin{align*}
            & \mathbb{E}[F_{x;t}(\omega_{0},\dots,\omega_{n})] \\
            & =(t/2)^{ -\frac{n}{2}+\frac{1}{2}\sum^n_{j=0}\mathrm{deg}(\omega'_j) }\Str_{x} \Big(\sum_{m = 1}^{n} (-2)^{m} \sum_{I \in \mathcal{P}_{m,n}}\mathbf{c}_{x}(\omega_{0}')\Phi^{D^2/2}_{t}\big(P(\omega_{I_1}),\dots, P(\omega_{I_m})(x,x)\big)\Big).
	\end{align*}
	As explained in the introduction, this identity plays a fundamental role in the geometry of the loop spaces.

	We now start with the preparation of the proof of Theorem \ref{local}. To this end, we fix once for all
	\begin{itemize}
		\item $d\coloneqq\dim(X)$,

		\item $l\coloneqq d/2$,

		\item $n\in\IN_{\geq 0}$,

		\item $\omega_{0}, \ldots, \omega_{n} \in \Omega_{\mathbb{T}}(X)$,

		\item $m_{n}\coloneqq\sum^{n}_{j=0}\mathrm{deg}(\omega'_{j})-n$,

		\item an oriented and smooth orthonormal frame $(e_{1},\dots, e_{d})$ for the tangent bundle of $X$, which is defined on some chart
			$U\subset X$,

		\item the induced coframe $(e^{1}, \ldots, e^{d})$ for the contangent bundle of $X$.
	\end{itemize}

	\begin{notation}
		1. If we are given $b\in [1,\infty]$, a family of random variables $f_{y}$, $y\in Y$, and for each $t>0$ a family of
		random variables $f_{y;t}$, $y\in Y$ (all taking values in a finite-dimensional normed space), we will say that
		\begin{equation*}
			\lim_{t\to 0+}f_{y;t}=f_{y}\quad\text{in $L^{b}_{y\in Y}(\IP)$,}
		\end{equation*}
		if the convergence is in $L^{b}(\IP)$, uniformly in $y\in Y$; in other words
		\begin{equation*}
			\lim_{t\to 0+}\sup_{y\in Y}\|f_{y;t}-f_{y}\|_{L^b(\IP)}=0.
		\end{equation*}
		2. Every matrix $A = (a_{ij}) \in \mathrm{Mat}(\IR;d\times d)$ induces an element
		$T(A)\in\Gamma_{C^\infty}(U,\End(\Sigma))$ via
		\begin{equation*}
			T_{x}(A) = \frac{1}{4} \sum_{1 \leq i, j \leq d}a_{ij}\mathbf{c}_{x}(e^{i}) \mathbf{c}_{x}(e^{j}),
		\end{equation*}
		and an element $\alpha(A)\in\Omega(X)$ via
		\begin{equation*}
			\alpha_{x}(A) \coloneqq \frac{1}{2} \sum_{i,j = 1}^{d} a_{ij}e^{i}_{x} \wedge e^{j}_{x}.
		\end{equation*}
		3. The symbol $R^{TX}$ denotes the curvature tensor of the Levi-Civita connection $\nabla^{TX}$ on the tangent bundle of $X$, and
		\begin{equation*}
			\Omega^{TX} = \frac{1}{2} \sum_{i,j = 1}^{d} R^{TX}(e_{i}, e_{j}) e^{i} \wedge e^{j}
		\end{equation*}
		denotes the curvature $2$-form of $R^{TX}$.
	\end{notation}

	\begin{remark}\label{aescd}
            Using the simple fact that
		\begin{equation*}
			\lim_{t \to 0+}t^{l}p(t,x,x) \mathrm{e}^{-\frac{1}{8}\int_0^t \mathrm{scal}(\mathsf{B}_s^{x;t}) \, \Id s}= (2 \pi)^{-l}
			\quad\text{in $L^{\infty}_{x\in X}(\IP)$},
		\end{equation*}
            we have
            \begin{align*}
                & \lim_{t\to 0+}F_{x;t}(\omega_{0},\dots,\omega_{n}) \\
                & = (2 \pi)^{-l}2^{-m_n/2}\lim_{t \to 0+}t^{\frac{m_{n}}{2} - l}\Str_{x} \Big(\sum_{m = 1}^{n} (-2)^{m} \sum_{I \in \mathcal{P}_{m,n}}\mathbf{c}(\omega_{0}') \int_{t\sigma_m}\delta\Psi^{x;t}_{\omega_{I_1}}(s_{1})\cdots \delta\Psi^{x;t}_{\omega_{I_m}}(s_{m}) \pa^{x;t}(t)^{-1}\Big),
		\end{align*}
            \begin{equation*}
                \text{in}\quad  \bigcap_{b\in [1,\infty)}L^{b}_{x\in X}(\IP),\quad\text{whenever the latter limit exists.}
            \end{equation*}
	\end{remark}

	We now determine the $t$-behaviour of the iterated Itô integrals in question. By definition, for every $I\in \mathcal{P}
	_{m,n}$, $1 \leq j \leq m$ we have a decomposition
	\begin{equation*}
		\Psi^{x;t}_{\omega_{I_j}}= \Psi^{x;t; 0}_{\omega_{I_j}}+ \Psi^{x,;t; 1}_{\omega_{I_j}},
	\end{equation*}
	where the first summand is an Itô integral and the second an ordinary Riemann integral.

	\begin{lemma}
		\label{It\^{o}} For all $b\geq 2$, $m\leq n$, $I\in \mathcal{P}_{m,n}$, and all
		$\nu=(\nu_{1},\dots,\nu_{m})\in \{0,1\}^{m}$, there exists a constant $C>0$ such that for all $0<t<1$,
		\begin{equation*}
			\sup_{x\in X}\mathbb{E}\left[\left \vert \int_{t \sigma_m}\delta \Psi^{x;t; \nu_1}_{\omega_{I_1}}(s_{1}) \cdots \delta
			\Psi^{x;t; \nu_m}_{\omega_{I_m}}(s_{m}) \right \vert^{b} \right] \leq C t^{\frac{b}{2}(m + |\nu|)},
		\end{equation*}
		where $|\nu|\coloneqq \nu_{1}+\cdots+\nu_{m}$.
	\end{lemma}

	\begin{proof}
		We are going to prove a slightly stronger statement, namely, the existence of $C>0$ such that for all $x\in X$,
		$0<t<1$ and $0<r\leq t$,
		\begin{equation}
			\label{sya}\mathbb{E}\left[\left \vert \int_{r \sigma_m}\delta \Psi^{x;t; \nu_1}_{\omega_{I_1}}(s_{1}) \cdots \delta
			\Psi^{x;t; \nu_m}_{\omega_{I_m}}(s_{m}) \right \vert^{b} \right] \leq C r^{\frac{b}{2}(m + \lvert \nu \rvert)}.
		\end{equation}
		Note that
		\begin{align*}
                & \Psi^{x;t; 0}_{\omega_{I_j}}(r)=\int^{r}_{0}\pa^{x;t}(s)^{-1}\mathrm{Sym}^{(1)}(P(\omega_{I_j}) )^{\flat}(\delta\mathsf{B}^{x;t}_{s})\pa^{x;t}(s), \\
                & \Psi^{x;t; 1}_{\omega_{I_j}}(r)=\int^{r}_{0}\pa^{x;t}(s)^{-1}V_{P(\omega_{I_j})}(\mathsf{B}^{x;t}_{s})\pa^{x;t}(s) \, \Id s.
		\end{align*}
		We proceed by induction over $m$: assume $m=1$. If $\nu_{1}=0$, we have
		\begin{align*}
                \mathbb{E}\left[\left \vert \int_{r \sigma_1}\delta \Psi^{x;t; \nu_1}_{\omega_{I_1}}(s_{1}) \right \vert^{b} \right] & =\mathbb{E}\left[\left \vert \int^{r}_{0}\pa^{x;t}(s)^{-1}\mathrm{Sym}^{(1)}(P(\omega_{I_1}) )^{\flat}(\delta\mathsf{B}^{x;t}_{s})\pa^{x;t}(s) \right \vert^{b}\right] \\
                & \leq C\; \mathbb{E}\left[ \left(\int^{r}_{0} \left \vert\mathrm{Sym}^{(1)}_{\mathsf{B}^{x;t}_s}(P(\omega_{I_j}) )^{\flat} \right \vert^{2} \, \Id s\right)^{\frac{b}{2}}\right] \\
                & \leq C\|\mathrm{Sym}^{(1)}(P(\omega_{I_1}) )^{\flat}\|_{\infty}^{b} r^{\frac{b}{2}},
		\end{align*}
		where we have used the Burkholder-Davis-Gundy inequality. If $\nu_{1}=1$, then
		\begin{align*}
                \mathbb{E}\left[\left \vert \int_{r \sigma_1}\delta \Psi^{x;t; \nu_1}_{\omega_{I_1}}(s_{1}) \right \vert^{b} \right] & =\mathbb{E}\left[\left \vert \int^{r}_{0}\pa^{x;t}(s)^{-1}V_{P(\omega_{I_1})}(\mathsf{B}^{x;t}_{s})\pa^{x;t}(s) \, \Id s \right \vert^{b}\right] \\
                & \leq \mathbb{E}\left[ \left(\int^{r}_{0}|V_{P(\omega_{I_1})}(\mathsf{B}^{x;t}_{s})| \, \Id s \right)^{b}\right] \\
                & \leq \|V_{P(\omega_{I_1})}\|_{\infty}^{b} r^{b}.
		\end{align*}
		Assume now (\ref{sya}) is correct for $m$. We have
		\begin{align}\label{po}
                & \mathbb{E}\left[\left \vert \int_{r \sigma_{m+1}}\delta \Psi^{x;t; \nu_1}_{\omega_{I_1}}(s_{1}) \cdots \delta \Psi^{x;t; \nu_m}_{\omega_{I_{m+1}}}(s_{m+1}) \right \vert^{b} \right] \\ \nonumber 
                & =\mathbb{E}\left[\left \vert \int^{r}_{0}\left(\int_{s \sigma_{m}}\delta \Psi^{x;t; \nu_1}_{\omega_{I_1}}(r_{1}) \cdots \delta \Psi^{x;t; \nu_m}_{\omega_{I_{m}}}(r_{m})\right)\delta \Psi^{x;t;\nu_{m+1}}_{\omega_{I_{m+1}}}(s) \right \vert^{b} \right].
		\end{align}
		If $\nu_{m+1}=0$, then (\ref{po}) is
		\begin{align*}
                & =\mathbb{E}\left[\left \vert \int^{r}_{0}\left(\int_{s \sigma_{m}}\delta \Psi^{x;t; \nu_1}_{\omega_{I_1}}(r_{1}) \cdots \delta \Psi^{x;t; \nu_m}_{\omega_{I_{m}}}(r_{m})\right)\pa^{x;t}(s)^{-1}\mathrm{Sym}^{(1)}(P(\omega_{I_{m+1}}) )^{\flat}(\delta\mathsf{B}^{x;t}_{s})\pa^{x;t}(s) \right \vert^{b} \right] \\
                & \leq C\mathbb{E}\left[\left( \int^{r}_{0}\left|\int_{s \sigma_{m}}\delta \Psi^{x;t; \nu_1}_{\omega_{I_1}}(r_{1}) \cdots \delta \Psi^{x;t; \nu_m}_{\omega_{I_{m}}}(r_{m})\right|^{2}|\mathrm{Sym}^{(1)}_{\mathsf{B}^{x;t}_s}(P(\omega_{I_{m+1}}) )^{\flat}|^{2} \, \Id s\right)^{\frac{b}{2}}\right] \\
                & \leq C\|\mathrm{Sym}^{(1)}(P(\omega_{I_{m+1}}) )^{\flat}\|_{\infty}^{b} \mathbb{E}\left[\left( \int^{r}_{0}\left|\int_{s \sigma_{m}}\delta \Psi^{x;t; \nu_1}_{\omega_{I_1}}(r_{1}) \cdots \delta \Psi^{x;t; \nu_m}_{\omega_{I_{m}}}(r_{m})\right|^{2} \, \Id s\right)^{\frac{b}{2}}\right] \\
                & =C\|\mathrm{Sym}^{(1)}(P(\omega_{I_{m+1}}) )^{\flat}\|_{\infty}^{b} r^{\frac{b}{2}}\mathbb{E}\left[\left( \int^{r}_{0}\left|\int_{s \sigma_{m}}\delta \Psi^{x;t; \nu_1}_{\omega_{I_1}}(r_{1}) \cdots \delta \Psi^{x;t; \nu_m}_{\omega_{I_{m}}}(r_{m})\right|^{2}\frac{\Id s}{r}\right)^{\frac{b}{2}}\right] \\
                & \leq C\|\mathrm{Sym}^{(1)}(P(\omega_{I_{m+1}}) )^{\flat}\|_{\infty}^{b} r^{\frac{b}{2}-1}\mathbb{E}\left[ \int^{r}_{0}\left|\int_{s \sigma_{m}}\delta \Psi^{x;t; \nu_1}_{\omega_{I_1}}(r_{1}) \cdots \delta \Psi^{x;t; \nu_m}_{\omega_{I_{m}}}(r_{m})\right|^{b} \, \Id s\right] \\
                & =C\|\mathrm{Sym}^{(1)}(P(\omega_{I_{m+1}}) )^{\flat}\|_{\infty}^{b} r^{\frac{b}{2}-1}\int^{r}_{0}\mathbb{E}\left[ \left|\int_{s \sigma_{m}}\delta \Psi^{x;t; \nu_1}_{\omega_{I_1}}(r_{1}) \cdots \delta \Psi^{x;t; \nu_m}_{\omega_{I_{m}}}(r_{m})\right|^{b}\right]\, \Id s \\
                & = C''\|\mathrm{Sym}^{(1)}(P(\omega_{I_{m+1}}) )^{\flat}\|_{\infty}^{b} r^{\frac{b}{2}(m +1+ \nu_1+\cdots+\nu_m+\nu_{m+1})},
		\end{align*}
		where we have used the Burkolder-Davis-Gundy inequality, Jensen's inequality and (\ref{sya}). If $\nu_{m+1}=1$, then
		(\ref{po}) is
		\begin{align*}
                & =\mathbb{E}\left[\left \vert \int^{r}_{0}\left(\int_{s \sigma_{m}}\delta \Psi^{x;t; \nu_1}_{\omega_{I_1}}(r_{1}) \cdots \delta \Psi^{x;t; \nu_m}_{\omega_{I_{m}}}(r_{m})\right)\pa^{x;t}(s)^{-1}V_{ P(\omega_{I_{m+1}}) }(\mathsf{B}^{x;t}_{s}) \pa^{x;t}(s) \, \Id s\right\vert^{b} \right] \\
                & \leq\|V_{P(\omega_{I_{m+1}})}\|_{\infty}^{b} \mathbb{E}\left[ \left(\int^{r}_{0}\left|\int_{s \sigma_{m}}\delta \Psi^{x;t; \nu_1}_{\omega_{I_1}}(r_{1}) \cdots \delta \Psi^{x;t; \nu_m}_{\omega_{I_{m}}}(r_{m})\right| \, \Id s\right)^{b} \right] \\
                & \leq\|V_{P(\omega_{I_{m+1}})}\|_{\infty}^{b} r^{b-1}\mathbb{E}\left[ \int^{r}_{0}\left|\int_{s \sigma_{m}}\delta \Psi^{x;t; \nu_1}_{\omega_{I_1}}(r_{1}) \cdots \delta \Psi^{x;t; \nu_m}_{\omega_{I_{m}}}(r_{m})\right|^{b} \, \Id s \right] \\
                & =\|V_{P(\omega_{I_{m+1}})}\|_{\infty}^{b} r^{b-1}\int^{r}_{0}\mathbb{E}\left[\left|\int_{s \sigma_{m}}\delta \Psi^{x;t; \nu_1}_{\omega_{I_1}}(r_{1}) \cdots \delta \Psi^{x;t; \nu_m}_{\omega_{I_{m}}}(r_{m})\right|^{b} \right] \, \Id s \\
                & \leq C\|V_{P(\omega_{I_{m+1}})}\|_{\infty}^{b} r^{b-1}\int^{r}_{0}s^{\frac{b}{2}(m + \nu_1+\cdots+\nu_m)}\, \Id s \\
                & = C'\|V_{P(\omega_{I_{m+1}})}\|_{\infty}^{b} r^{b+\frac{b}{2}(m + \nu_1+\cdots+\nu_m)} \\
                & = C'\|V_{P(\omega_{I_{m+1}})}\|_{\infty}^{b} r^{\frac{b}{2}(m+1 + \nu_1+\cdots+\nu_m+\nu_{m+1})},
		\end{align*}
		where we have used Jensen's inequality. This completes the proof.
	\end{proof}

	\begin{remark}
		Assume $m_{n}/2 > l$. Then we have
		\begin{equation*}
			\deg(\omega_{0}') + \ldots + \deg(\omega_{n}') > d + n,
		\end{equation*}
		implying
		\begin{equation*}
			( \hat{A}(X)\wedge \omega_{0}' \wedge \omega_{1}'' \wedge \dots \wedge \omega_{n}'')^{[d]}= 0.
		\end{equation*}
		Moreover, the above estimates show
		\begin{equation*}
			\lim_{t\to 0+}F_{x,t}(\omega_{0}, \ldots, \omega_{n})= 0\quad\text{in
			$\bigcap_{b\in [1,\infty)}L^{b}_{x\in X}(\IP)$,}
		\end{equation*}
		and we see that all claims made in Theorem \ref{local} hold under this degree assumption.
	\end{remark}

	The following asymptotic behaviour of the stochastic parallel transport on the tangent bundle has been proven by Hsu in
	\cite{hsu}:

	\begin{proposition}
		Let $t\in [0,1]$ and let $\pa^{x;t}_{TX}$ denote the parallel transport along $\mathsf{B}
		^{x;t}$ with respect to $\nabla^{TX}$, and for every $x\in U$ define a process
		$\mathsf{C}^{x;t}(s)\in\mathrm{SO}(\IR;d)$, $s\in [0,t]$, by
		\begin{equation*}
			\mathsf{C}^{x;t}_{ij}(s)\coloneqq(\pa^{x;t}_{TX}(s)^{-1}e_{i}^{x},e_{j}^{x}).
		\end{equation*}
		Then for all $x\in U$ there exists a standard Euclidean Brownian bridge $\mathsf{Y}^{x}_{s}\in T_{x} X$, $s\in [0,1]$,
		such that for all $b\geq 1$ one has
		\begin{equation*}
			\lim_{t\to 0+}\frac{1}{t}(\mathsf{C}^{x;t}(s)-1)=-\frac{1}{2} \int_{0}^{1} \big( R^{TX}(e^{x}, e^{x})\mathsf{Y}^{x}_{r},
			\delta \mathsf{Y}^{x}_{r} \big) \quad\text{in $L^{b}_{(x,s)\in U\times [0,t]}(\IP)$},
		\end{equation*}
		where the random variable
		$\int_{0}^{1} \big( R^{TX}(e^{x}, e^{x})\mathsf{Y}^{x}_{r}, \delta \mathsf{Y}^{x}_{r} \big)\in \mathfrak{so}(\IR;d)$ is
		defined by
		\begin{equation*}
			\Big(\int_{0}^{1} \big( R^{TX}(e^{x}, e^{x})\mathsf{Y}^{x}_{s}, \delta \mathsf{Y}^{x}_{s} \big)\Big)_{ij}\coloneqq\int_{0}
			^{1} \big( R^{TX}(e^{x}_{i}, e^{x}_{j})\mathsf{Y}^{x}_{r}, \delta \mathsf{Y}^{x}_{r} \big).
		\end{equation*}
	\end{proposition}

	We record the following simple consequence of the above result:

	\begin{lemma}\label{Limit parallel transport} 
		\emph{a)} One has
		\begin{align}
			\label{asss}\lim_{t\to 0+}\sup_{(x,s)\in U\times [0,t]}\left|\mathsf{C}^{x;t}(s)-1\right|= 0\quad\text{a.s.}
		\end{align}
		\emph{b)} For all $b\geq 1$, $\delta>0$, there exists a constant $C>0$ such that for all $t\in [0,1]$,
		\begin{equation*}
			\sup_{(x,s)\in U\times [0,t]}\mathbb{P}\{\left|\mathsf{C}^{x;t}(s)-1\right| \geq \delta\}\leq Ct^{b}.
		\end{equation*}
		\emph{c)} Assume that for all $x\in U$, $t\in [0,1]$ one has
		\begin{equation*}
			\sup_{s\in [0,t]}\left|\mathsf{C}^{x;t}(s)-1\right|< 1\quad\text{ a.s.},
		\end{equation*}
		so that
		\begin{equation*}
			\mathsf{A}^{x;t}(s)\coloneqq\log(\mathsf{C}^{x;t}(s))=\sum^{\infty}_{k=1}\frac{(-1)^{k}}{k}(\mathsf{C}^{x;t}(s)-1)^{k}
			,\quad s\in [0,t],
		\end{equation*}
		is the uniquely determined process $\mathsf{A}^{x;t}(s) \in \mathfrak{so}(\IR;d)$, $s\in [0,t]$, with
		$\mathsf{C}^{x;t}(s) = \mathrm{e}^{\mathsf{A}^{x;t}(s)}$. Then for all $b\geq 1$,
		\begin{equation*}
			\lim_{t\to 0+}\frac{1}{t}\mathsf{A}^{x;t}(s)=-\frac{1}{2} \int_{0}^{1} \big( R^{TX}(e^{x}, e^{x})\mathsf{Y}^{x}_{r}, \delta
			\mathsf{Y}^{x}_{r}) \quad\text{in $L^{b}_{(x,s)\in U\times [0,t]}(\IP)$,}
		\end{equation*}
		in particular,
		\begin{equation*}
			\lim_{t\to 0+}\frac{1}{t}T_{x}(\mathsf{A}^{x;t}(s))=-\sum^{d}_{i,j = 1}\frac{1}{8} \int_{0}^{1} \big( R^{TX}(e^{x}_{i},
			e^{x}_{j})\mathsf{Y}^{x}_{r}, \delta \mathsf{Y}^{x}_{r} \big) \mathbf{c}(e^{x}_{i})\mathbf{c}(e^{x}_{j})\quad\text{in
			$L^{b}_{(x,s)\in U\times [0,t]}(\IP)$.}
		\end{equation*}

		\emph{d)} In the situation of \emph{c)}, for all $x\in U$, $t\in [0,1]$, $s\in[0,t]$, one has
		$\pa^{x;t}(s)=\mathrm{e}^{T_x(\mathsf{A}^{x;t}(s))}$ a.s.
	\end{lemma}

	\begin{proof}
		a) The convergence from the previous proposition implies the existence of $C'>0$, depending on $b$, such that
		\begin{align}
			\label{evb}\sup_{(x,s)\in U\times [0,t]}\mathbb{E}\left[\left|\mathsf{C}^{x;t}(s)-1\right|^{b}\right]\leq C' t^{b},\quad\text{in particular,}\quad \lim_{t\to 0+}\sup_{(x,s)\in U\times [0,t]}\mathbb{E}\left[\left|\mathsf{C}^{x;t}(s)-1\right|\right]=0.
		\end{align}
		This observation easily implies (\ref{asss}): indeed, assuming the contrary, we find a sequence of times $t_{n}\to 0$
		such that on every set of full measure, the sequence of random variables
		\begin{equation*}
			\sup_{(x,s)\in U\times [0,t]}\left|\mathsf{C}^{x;t_n}(s)-1\right|\quad\text{ does not converge to $0$}.
		\end{equation*}
		Thus, there exists $\epsilon>0$ and a subsequence $t_{n_k}$ of $t_{n}$ such that for all $k$ one has
		\begin{equation*}
			\sup_{(x,s)\in U\times [0,t]}\left|\mathsf{C}^{x;t_{n_k}}(s)-1\right|\geq \epsilon,
		\end{equation*}
		which contradicts the limit in (\ref{evb}).

		b) By Markov's inequality and the first estimate in (\ref{evb}),
		\begin{equation*}
			\mathbb{P}\{\left|\mathsf{C}^{x;t}(s)-1\right| \geq \delta\}=\mathbb{P}\{\left|\mathsf{C}^{x}(t)-1\right|^{b} \geq
			\delta^{b}\}\leq \mathbb{E}\left[\left|\mathsf{C}^{x;t}(s)-1\right|^{b}\right] /\delta^{b}\leq (C't^{b})/\delta^{b}
			,
		\end{equation*}
		where $C'$ is uniform in $(x,s)\in U\times [0,t]$.

		c) For all $s\in [0,t]$ one a.s. has
		\begin{align*}
                & \frac{1}{t}\mathsf{A}^{x;t}(s)+\frac{1}{2} \int_{0}^{1} \big( R^{TX}(e^{x}, e^{x})\mathsf{Y}^{x}_{r}, \delta \mathsf{Y}^{x}_{r} \big) \\
                & =\frac{1}{t}\sum^{\infty}_{k=1}\frac{(-1)^{k}}{k}(\mathsf{C}^{x;t}(s)-1)^{k} +\frac{1}{2} \int_{0}^{1} \big( R^{TX}(e^{x}, e^{x})\mathsf{Y}^{x}_{r}, \delta \mathsf{Y}^{x}_{r} \big) \\
                & =\frac{1}{t}\sum^{\infty}_{k=2}\frac{(-1)^{k}}{k}(\mathsf{C}^{x;t}(s)-1)^{k} +\frac{1}{t}(\mathsf{C}^{x;t}(s)-1)+\frac{1}{2} \int_{0}^{1} \big( R^{TX}(e^{x}, e^{x})\mathsf{Y}^{x}_{r}, \delta \mathsf{Y}^{x}_{r} \big),
		\end{align*}
		and so
		\begin{align*}
                & \sup_{(x,s)\in U\times [0,t]}\mathbb{E}\Big[\Big|\frac{1}{t}\mathsf{A}^{x;t}(s)+\frac{1}{2} \int_{0}^{1} \big( R(e^{x}, e^{x})\mathsf{Y}^{x}_{r}, \delta \mathsf{Y}^{x}_{r} \big)\Big|^{b}\Big] \\
                & \leq 2^{b-1}\sup_{(x,s)\in U\times [0,t]}\mathbb{E}\Big[\Big|\sum^{\infty}_{k=2}\frac{1}{t}\frac{(-1)^{k}}{k}(\mathsf{C}^{x;t}(s)-1)^{k}\Big|^{b}\Big] \\
                & \quad+2^{b-1}\sup_{(x,s)\in U\times [0,t]}\mathbb{E}\Big[\Big|\frac{1}{t}(\mathsf{C}^{x;t}(s)-1)+\frac{1}{2} \int_{0}^{1} \big( R^{TX}(e^{x}, e^{x})\mathsf{Y}^{x}_{r}, \delta \mathsf{Y}^{x}_{r} \big) \Big|^{b}\Big].
		\end{align*}

		The second term converges to $0$ as $t\to 0$+; using the Minkowski inequality and the first inequality in (\ref{evb}),
		the first term can be controlled as follows,
		\begin{align*}
                & \sup_{(x,s)\in U\times [0,t]}\mathbb{E}\Big[\Big|\sum^{\infty}_{k=2}\frac{1}{t}\frac{(-1)^{k}}{k}(\mathsf{C}^{x;t}(s)-1)^{k}\Big|^{b}\Big] \\
                & \qquad \leq \Big(\sum^{\infty}_{k=2}\Big(\sup_{(x,s)\in U\times [0,t]}\mathbb{E}\Big[\Big|\frac{1}{t}\frac{(-1)^{k}}{k}(\mathsf{C}^{x;t}(s)-1)^{k}\Big|^{b}\Big]\Big)^{\frac{1}{b}}\Big)^{b} \\
                & \qquad \leq C \Big( \sum^{\infty}_{k=1}t^{k}\Big)^{b} \to 0,
		\end{align*}
		as $t\to 0+$.

		d) This follows easily from the defining relation of stochastic parallel transport (\ref{para}) and the fact that if
		locally $\nabla^{TX}$ is given by
		\begin{equation*}
			\nabla^{TX}_{\partial_i}e_{j}=\sum_{k}\omega^{k}_{ij}e_{k},\quad \omega^{k}_{ij}\in C^{\infty}(U),
		\end{equation*}
		then the spin Levi-Civita connection is given by
		\begin{equation*}
			\nabla_{\partial_i}=\partial_{i}+\frac{1}{4}\sum_{jk}\sum_{k}\omega^{k}_{ij}\mathbf{c}(e_{j})\mathbf{c}(e_{k}).
		\end{equation*}
		This completes the proof.
	\end{proof}

	\begin{remark}\label{lplp}
		1. In view of this result, writing
		\begin{align*}
                & t^{\frac{m_{n}}{2} - l}\Str_{x} \Big(\sum_{m = 1}^{n} (-2)^{m} \sum_{I \in \mathcal{P}_{m,n}}\mathbf{c}(\omega_{0}') \int_{t\sigma_m}\delta\Psi^{x;t}_{\omega_{I_1}}(s_{1})\cdots \delta\Psi^{x;t}_{\omega_{I_m}}(s_{m}) \pa^{x;t}(t)^{-1}\Big) \\
                & = t^{\frac{m_{n}}{2} - l}1_{\{\sup_{s\in [0,t]}\left|\mathsf{C}^{x;t}(s)-1\right| \geq 1\}}\Str_{x} \Big(\sum_{m = 1}^{n} (-2)^{m} \sum_{I \in \mathcal{P}_{m,n}}\mathbf{c}(\omega_{0}') \times \\ 
                & ~~~~~~~~~~~~~~~~~~~~~~~~~~~~~~~~~~ \times \int_{t\sigma_m}\delta\Psi^{x;t}_{\omega_{I_1}}(s_{1})\cdots \delta\Psi^{x;t}_{\omega_{I_m}}(s_{m}) \pa^{x;t}(t)^{-1}\Big) \\
                & +t^{\frac{m_{n}}{2} - l}1_{\{\sup_{s\in [0,t]}\left|\mathsf{C}^{x;t}(s)-1\right| < 1\}}\Str_{x} \Big(\sum_{m = 1}^{n} (-2)^{m}\sum_{I \in \mathcal{P}_{m,n}}\mathbf{c}(\omega_{0}') \times \\ 
                & ~~~~~~~~~~~~~~~~~~~~~~~~~~~~~~~~~~ \times \int_{t\sigma_m}\delta\Psi^{x;t}_{\omega_{I_1}}(s_{1})\cdots \delta\Psi^{x;t}_{\omega_{I_m}}(s_{m}) \pa^{x;t}(t)^{-1}\Big),
		\end{align*}
		and using Lemma \ref{It\^{o}}, Lemma \ref{Limit parallel transport}, and Remark \ref{aescd} shows that for all
		$t\in (0,1]$, $x\in U$ there exists a uniquely determined random variable $\mathsf{A}^{x}(t) \in \mathfrak{so}(\IR;d)$, with
		\begin{equation*}
                \pa^{x;t}(t)^{-1}= \mathrm{e}^{T_x(\mathsf{A}^x(t))},
		\end{equation*}
		and 
		\begin{align*}
                &\lim_{t\to 0+}F_{x;t}(\omega_{0},\dots,\omega_{n}) \\
                & = (2 \pi)^{-l}2^{-m_n/2}\lim_{t\to 0+}t^{\frac{m_{n}}{2} - l}\Str_{x} \Big(\sum_{m = 1}^{n} (-2)^{m} \sum_{I \in \mathcal{P}_{m,n}}\mathbf{c}(\omega_{0}') \int_{t\sigma_m}\delta\Psi^{x;t}_{\omega_{I_1}}
			(s_{1})\cdots \delta\Psi^{x;t}_{\omega_{I_m}}(s_{m}) \mathrm{e}^{T_x(\mathsf{A}^x(t))}\Big),
		\end{align*}
		in
            \begin{equation*}
                \bigcap_{b\in [1,\infty)}L^{b}_{x\in U}(\IP),\quad\text{whenever the limit on the right hand side exists.}
            \end{equation*}
            2. In view of Lemma \ref{Limit parallel transport} and the Taylor expansion
		\begin{align*}
			\pa^{x;t}(t)^{-1}= \sum_{k = 0}^{l}\frac{T_x(\mathsf{A}^x(t))^{k}}{k!} + \mathsf{R}^{x}_l(t),
		\end{align*}
                for all $i\in\IN$ there exists $C_{l,i}>0$ such that for all $0<t<1$ the remainder satisfies
		\begin{align}\label{Estimate remainder}
                \sup_{x\in U}\mathbb{E}[\lvert \mathsf{R}^{x}_l \rvert^{i}] \leq C_{l,i}t^{i(l + 1)},
		\end{align}
		and for all $i\in\IN$ there exists $C_i>0$ such that for all $0<t<1$,
		\begin{equation*}
		  \sup_{x\in U}	\mathbb{E}[\lvert T_x(\mathsf{A}^x(t))^{i} \rvert]\leq C_{i} t^{i}.
		\end{equation*}
		Note that $m_{n}/2 < 0$ only if one of the forms $\omega_{j}'$, $1 \leq j \leq n$, has degree $0$. However, for such
		$\omega_{j}'$ the operators $P(\omega_{j}')$, $P(\omega_{j-1}', \omega_{j}')$ and $P(\omega_{j}', \omega_{j+1}')$ will
		be zero. From this and \eqref{Estimate remainder}, we deduce
		\begin{align*}\lim_{t\to 0+}F_{x;t}(\omega_{0},\dots,\omega_{n}) & = (2 \pi)^{-l}2^{-m_n/2}\lim_{t \to 0+}t^{m_n/2 - l}\Str\Bigg(\sum_{m = 1}^{n} (-2)^{m} \sum_{I \in \mathcal{P}_{m,n}}\mathbf{c}_{x}(\omega_{0}') \\
                & \quad \vphantom{\sum_{m = 1}^n}\mathbb{E}\left[ \int_{t\sigma_m}\delta\Psi^{x;t}_{\omega_{I_1}}(s_{1})\cdots \delta\Psi^{x;t}_{\omega_{I_m}}(s_{m}) \sum_{k = 0}^{l}\frac{T_x(\mathsf{A}^x(t))^{k}}{k!}\right]\Bigg)\quad\text{in $\bigcap_{b\in [1,\infty)}L^b_{x\in U}(\IP)$},
            \end{align*}
            whenever the latter limit exists.

            3. The above estimates show that the only summands in the above formula contributing to the limit are those for which
		\begin{equation}\label{Fundamental inequality order}
                \frac{m_{n}}{2}+ k + \frac{m + \lvert \nu \rvert}{2}\leq l,\quad \text{where $\lvert \nu \rvert = \nu_{1} + \dots + \nu_{m}$.}
		\end{equation}
	\end{remark}

	An algebraic argument will show that all summands for which $\frac{m_{n}}{2}+ k + \frac{m + \lvert \nu \rvert}{2}< l$ 
        vanish already at the level of loops (so prior to taking the limit and expectation) and even more, 
        that we only get contributions for $\nu = (1, \ldots, 1) \in \{0,1\}^{n}$, i. e. when $m = n$. This will be implemented by the following algebraic machinery:

	\begin{definition}
		\label{patodifilt} Given $x\in U$, $r \in \IN_{\geq 1}$ we say that $T\in \End(\Sigma_{x})$ has \emph{Patodi order} $\leq
		r$, if there exist $1 \leq i \leq \lfloor r \rfloor$ and matrices $A_{1}, \ldots, A_{i}\in \mathfrak{so}(\IR;d)$
		such that $T = T_{x}(A_{1}) \cdots T_{x}(A_{i})$. We simply write $\mathrm{Pat}(T) \leq r$ in this situation.
	\end{definition}

	\begin{remark}
		For all $x\in U$ the Patodi order defines a filtration of $\End(\Sigma_{x})$ via
		\begin{equation*}
			\End^{r}(\Sigma_{x})\coloneqq\{T\in\End(\Sigma_{x})\colon \mathrm{Pat}(T) \leq r\}.
		\end{equation*}
	\end{remark}

	This notion is justified by the following version of Patodi's lemma, which we have borrowed from \cite{hsu} (cf. Lemma
	7.4.3 therein).

	\begin{lemma}
		Let $x\in U$.
		
		\emph{1.} For every $T \in \End(\Sigma_{x})$ with $\mathrm{Pat}(T) \leq l-1$ one has $\Str_{x}(T) = 0$.
		
		\emph{2.} For every collection of $l$ many matrices $A_{1}, \ldots, A_{l} \in \mathfrak{so}(\IR;d)$ one has
		\begin{equation*}
			\Str_{x}\big(T_{x}(A_{1}) \cdots T_{x}(A_{l})\big) e^{1}_{x} \wedge \dots \wedge e^{d}_{x} = \frac{1}{\sqrt{-1}^{l}}
			\alpha_{x}(A_{1}) \wedge \dots \wedge \alpha_{x}(A_{l}).
		\end{equation*}
	\end{lemma}

	\begin{proof}[Proof of Theorem \ref{local}]
            We fix $b\in [1,\infty)$ and understand all limits in the sequel with respect to $L^b_{x\in U}(\IP)$. We may assume that the forms $\omega_{0}, \ldots, \omega_{n}$ are of pure degree. Since $\hat{A}(X)$ only
		has non-vanishing components in degrees which are multiples of $4$ and since $d$ is even, we have
		\begin{equation*}
			\big(\hat{A}(X) \wedge \omega_{0}' \wedge \omega_{1}'' \wedge \dots \wedge \omega_{n}''\big)^{[d]}_{x} = 0,\>\>\text{if $\deg(\omega_{0}') + \deg(\omega_{1}'') + \dots + \deg(\omega_{n}'')$ is odd,}
		\end{equation*}
            which happens exactly when an odd
		number of those forms are odd. One checks that in this case, $m_{n}$ is odd which implies that the endomorphisms 
            from the definition of $F_{x;t}(\omega_0,\dots,\omega_n)$ are odd so that their supertrace vanishes. 
            Therefore, we may assume that an even number of the forms $\omega_{0}', \omega_{1}'', \ldots, \omega_{n}''$ is odd.

		1. We start with the case $n = 0$. Let $\omega\coloneqq\omega_{0}$ and let $\deg(\omega')$ be the even number $r = 2q$. In view of Remark \ref{lplp} we have
            \begin{equation*}
                \lim_{t\to 0+}F_{x;t}(\omega)=(
			2 \pi)^{-l}2^{-q}\lim_{t \to 0+}t^{q-l}\Str_x\Bigg(\mathbf{c}_{x}(\omega') 
			\frac{T_x(\mathsf{A}^x(t))^{l-q}}{(l-q)!}\Bigg),
            \end{equation*}
            in case of existence of the latter, which in turn is implied by Lemma \ref{Limit parallel transport}, giving 
		\begin{equation*}
                F_x(\omega) = \lim_{t\to 0+}F_{x;t}(\omega) = (
			2 \pi)^{-l}2^{-q}\Str_x\Bigg(\mathbf{c}_{x}(\omega') 
			\frac{T_x(\mathsf{D}^x)^{l-q}}{(l-q)!}\Bigg),
		\end{equation*}
            where
            \begin{equation*}
			\mathsf{D}_{x} \coloneqq -\frac{1}{2} \int_{0}^{1} \big( R^{TX} (e_{i}^{x}, e_{j}^{x})\mathsf{Y}_{s}^{x}, \delta \mathsf{Y}_{s}^{x}
			\big).
		\end{equation*}
            It remains to evaluate the expectation of $F_x(\omega) e^1_x\wedge \cdots \wedge e_x^d$. To this end, we expand
		\begin{equation*}
			\omega'_{x}= \sum_{1 \leq i_1 < \dots < i_r \leq d}\alpha_{i_1, \ldots, i_r}e^{i_1}_{x} \wedge \dots \wedge e^{i_r}
			_{x},
		\end{equation*}
		for certain $\alpha_{i_1, \ldots, i_r}\in \mathbb{R}$. One has
		\begin{equation*}
			2T_x(A^{ij})=\mathbf{c}_{x}(e^{i})\mathbf{c}_{x}(e^{j}) \in \End(\Sigma_{x}),
		\end{equation*}
		where $A^{ij}\in \mathfrak{so}(\IR;d)$ is given by $a_{ij}= 1$, $a_{ji}= -1$ with all other entries equal to $0$. More
		generally, if $\tau$ is a permutation of $\{1, \dots, r\}$ and $1 \leq i_{1} < \dots < i_{r} \leq d$, then
		\begin{equation*}
			\mathbf{c}_{x}\big(e^{i_{\tau(1)}}\big) \cdots \mathbf{c}_{x}\big(e^{i_{\tau(r)}}\big) = 2^{q} T_x(A^{\tau(1) \tau(2)})
			\cdots T_x(A^{\tau(r-1) \tau(r)}).
		\end{equation*}
		Using this and the fact that
		\begin{align*}
                \mathbf{c}_{x}\big(\alpha_{i_1, \ldots, i_r}e^{i_1}\wedge \dots \wedge e^{i_r}\big) & = \frac{\alpha_{i_1, \ldots, i_r}}{r!}\sum_{\tau \in S_r}\sgn(\tau) \mathbf{c}_{x}(e^{i_{\tau(1)}}) \cdots \mathbf{c}_{x}(e^{i_{\tau(r)}})\\
			& = \alpha_{i_1, \ldots, i_r}\mathbf{c}_x(e^{i_1}) \cdots \mathbf{c}_x(e^{i_r}),
		\end{align*}
		we may assume $\omega'_{x}= e_{x}^{1} \wedge \dots \wedge e_{x}^{r}$. Let $A_{1}, \ldots, A_{q} \in \mathfrak{so}(\IR;d)$ be the matrices corresponding to
		$\mathbf{c}_{x}(e^{1}) \mathbf{c}_{x}(e^{2}), \ldots, \mathbf{c}_{x}(e^{r-1}) \mathbf{c}_{x}(e^{r})$ via
		$A\mapsto T_x(A)$.  
            With the help of the second part of Patodi's Lemma one can calculate
		\begin{align*}
                & \mathbb{E}\left[F_x(\omega)\right] \cdot \underline{\mu}_x \\
                & = \mathbb{E}\left[F_x(\omega) e^1_x\wedge \cdots \wedge e_x^d \right] \\
                & = (2 \pi)^{-l}2^{-q}\mathbb{E}\left[\Str\left(\mathbf{c}_{x}(\omega') \frac{T(\mathsf{D}^{x})^{l-q}}{{(l-q)!}}\right)e^{1}_{x} \wedge \dots \wedge e^{d}_{x}\right] \\
                & = \frac{1}{\big(2 \pi \sqrt{-1}\big)^{l}(l-q)!}\mathbb{E}\big[\alpha_x(A_{1}) \wedge \dots \wedge \alpha_x(A_{q}) \wedge \alpha_x(\mathsf{D}^{x})^{\wedge (l-q)}\big] \\
                & = \frac{1}{(2 \pi \sqrt{-1})^{l}}\mathbb{E}\Big[(\omega'_{x} \wedge \mathrm{e}^{\mathsf{J}_x})^{[d]}\Big] \\
                & = \frac{1}{(2 \pi \sqrt{-1})^{l}}(\omega'_{x} \wedge \mathbb{E}[\mathrm{e}^{\mathsf{J}_x}])^{[d]},\quad\text{where $\mathsf{J}_{x} \coloneqq - \int_{0}^{1} ( \Omega^{TX}_{x} \mathsf{Y}^{x}_{s}, \delta \mathsf{Y}^{x}_{s} )$.}
		\end{align*}
		A well-known variant of Lévy's stochastic area formula states $\mathbb{E}[\mathrm{e}^{\mathsf{J}_x}] = \hat{A}(X)_{x}$ 
            (see e. g. Lemma 7.6.6 in \cite{hsu}), and the proof of the $n=0$ case is complete.

		2. We proceed with $n = 1$. We may assume that both $\omega_{0}'$ and $\omega_{1}''$ are either even or odd so that in
		particular, $m_{1}$ is even. We have
		\begin{equation*}
			\lim_{t \to 0+}F_{x;t}(\omega_{0}, \omega_{1}) = -(2 \pi)^{-l}2^{-\frac{m_1}{2}  + 1}\lim
			_{t \to 0+}t^{\frac{m_1}{2} -l }\Str_x\left(\mathbf{c}_{x}(\omega_{0}')\Psi^{x;t}_{\omega_1}(t) \sum
			_{k = 0}^{l}\frac{T_x(\mathsf{A}^x(t))^{k}}{k!}\right),
		\end{equation*}
		where due to \eqref{Fundamental inequality order} and the fact that $\frac{m_{1}}{2}$ is an integer the only summands contributing are those for which
		\begin{equation*}
			\frac{m_{1}}{2}+ k \leq l-1.
		\end{equation*}
		A standard calculation (see e. g. \cite{boldt}) shows that
		\begin{align}\label{Principal symbol curvature}
                \begin{split}
                    \mathrm{Sym}^{(1)}(P(\omega_{1}))^{\flat}(X) & = - 2 \mathbf{c}(X \lrcorner \omega_{1}') \quad \text{for all $X \in \Gamma_{C^\infty}(X, TX)$, and} \\
                    \quad V_{P(\omega_1)} & = \mathbf{c}(\Id^{\dagger} \omega_{1}') + \mathbf{c}(\omega_{1}'').
                \end{split}
		\end{align}
		It follows that
		\begin{equation*}
			\mathrm{Pat}\big(\mathbf{c}_{x}(\omega_{0}') \Psi^{x;t;0}_{\omega_1}T_{\mathsf{A}^x(t)}^{k} \big) \leq \frac{\omega_{0}'}{2}
			+ \frac{\omega_{1}' - 1}{2}+ k = \frac{m_{1}}{2}+ k \leq l - 1
		\end{equation*}
		and hence, by the first part of Patodi's lemma,
		\begin{equation*}
                \lim_{t \to 0+}F_{x;t}(\omega_{0}, \omega_{1})  = -(2 \pi)^{-l}2^{-m_1/2 + 1}\lim
                _{t \to 0+}t^{m_1 / 2 -l }\Str\left(\mathbf{c}_{x}(\omega_{0}')\Psi^{x;t;1}_{\omega_1}(t) \sum
			_{k = 0}^{l^*}\frac{T_x(\mathsf{A}^x(t))^{k}}{k!}\right) 
            \end{equation*}
		where $l^{*} \coloneqq l - 1 - \frac{m_{1}}{2}$. Similarly, we see that
		\begin{equation*}
			\mathrm{Pat}\Big(\mathbf{c}_{x}(\omega_{0}') \int_{0}^{t} \pa^{x;t}(s)^{-1}\mathbf{c}_{\mathsf{B}_s^{x;t} }(\Id^{\dagger}
			\omega_{1}') \pa^{x;t}(s) \, \Id s \, T_x(\mathsf{A}^x(t))^{k}\Big) \leq \frac{m_{1}}{2}+ k \leq l - 1
		\end{equation*}
		for $0 \leq k \leq l^{*}$ and
		\begin{equation*}
			\mathrm{Pat}\Big(\mathbf{c}_{x}(\omega_{0}') \int_{0}^{t} \pa^{x;t}(s)^{-1}\mathbf{c}_{\mathsf{B}_s^{x;t} }(\omega_{1}
			'') \pa^{x;t}(s) \, \Id s \, T_x(\mathsf{A}^x(t))^{k} \Big) \leq \frac{\omega_{0}'}{2}+ \frac{\omega_{1}' + 1}{2}+ k
			= \frac{m_{1}}{2}+ 1 + k \leq l - 1
		\end{equation*}
		for $0 \leq k \leq l^{*} - 1$. Consequently,
		\begin{align*}
                & \lim_{t \to 0+}F_{x;t}(\omega_{0}, \omega_{1}) \\
                & = -(2 \pi)^{-l}2^{-\frac{m_1}{2}  + 1}\lim_{t \to 0+}\Str_x\left(\mathbf{c}_{x}(\omega_{0}')\frac{1}{t} \int_{0}^{t} \pa^{x;t}(s)^{-1}\mathbf{c}_{\mathsf{B}_s^{x;t} }( \omega_{1}'') \pa^{x;t}(s) \, \Id s\frac{T_x(\mathsf{A}^x(t))^{l^*}}{t^{l^*}l^{*}!}\right) \\
                & = -(2 \pi)^{-l}2^{-\frac{m_1}{2}  + 1}\Str_x\left(\mathbf{c}_{x}(\omega_{0}')\mathbf{c}_{x}( \omega_{1}'') \frac{T_x(\mathsf{D}^x)^{l^*}}{l^{*}!}\right),
		\end{align*}
  were the latter equality follows from Lemma \ref{Limit parallel transport} and the fact that the integral is an ordinary Riemann integral. Arguing in the same way as for $n = 0$ employing the second part of Patodi's lemma, from this it follows that
		\begin{align*}
                & \mathbb{E}[F_x(\omega_0, \omega_1)] \cdot \underline{\mu}_x \\
                & \qquad = -(2 \pi)^{-l}2^{-\frac{m_1}{2}  + 1}\mathbb{E}\left[\Str_x\left(\mathbf{c}_{x}(\omega_{0}')\mathbf{c}_{x}( \omega_{1}'') \frac{T_x(\mathsf{D}^x)^{l^*}}{l^{*}!}\right) e^1_x \wedge\cdots\wedge e^d_x \right] \\
                & \qquad = -\frac{2^{2}}{(2 \pi \sqrt{-1})^{l}}\big(\hat{A}(X) \wedge \omega_{0}' \wedge \omega_{1}''\big)_{x}^{[d]},
		\end{align*}
		as desired.

		3. Finally, we may deal with the general case $n \geq 2$.

		To this end, let $1 \leq m \leq n - 1$ and $\nu \in \{0,1\}^{m}$ and $I \in \mathcal{P}_{m,n}$. We want to estimate
		the Patodi orders of the summands in
		\begin{equation}
			\label{Patodi order to be estimated}\mathbf{c}_{x}(\omega_{0}') \int_{t \sigma_m}\delta \Psi^{x;t; \nu_1}_{\omega_{I_1}}
			(s_{1}) \cdots \delta \Psi^{x;t; \nu_m}_{\omega_{I_m}}(s_{m}) \, \frac{T_x(\mathsf{A}^x(t))^{k}}{k!}.
		\end{equation}
		First of all, we employ \eqref{Components P} to see that in order for the above integral to be non-zero we
		necessarily have that $\# I_{j} = 2$ for $n - m$ indices $j$ and $\# I_{j} = 1$ for $2m - n$ indices $j$. However,
		since $P(\cdot, \cdot)$ maps to zero-order differential operators, we have
		\begin{equation*}
			\mathrm{Sym}^{(1)}(P(\omega_{I_j}))^{\flat} = 0, \quad\text{ if $\# I_{j} = 2$}.
		\end{equation*}

		Thus, to get a non-zero integral above, we must have $\nu_{j} = 1$ precisely for those indices $j$. However, a
		standard calculation shows that for $\eta_{1}, \eta_{2} \in \Omega(X)$ of pure degree, the endomorphism
		\begin{equation*}
			\mathbf{c}_{x}(\eta_{1} \wedge \eta_{2}) - \mathbf{c}_{x}(\eta_{1}) \mathbf{c}_{x}(\eta_{2})\in\End(\Sigma_{x})
		\end{equation*}
		can be written as a linear combination of $\deg( \eta_{1}) + \deg( \eta_{2}) - 2$ many Clifford multiplications.
		Together with \eqref{Principal symbol curvature} this implies that for $0 \leq k \leq l$, the Patodi order of every
		summand in \eqref{Patodi order to be estimated} is less or equal to
		\begin{align*}
                & \dfrac{\deg(\omega_0') + \dots + \deg(\omega_n') -2(n-m) + (\lvert \nu \rvert - (n-m)) - (2m - n - (\lvert \nu \rvert - (n-m))}{2}+ k \\
                & \qquad = \dfrac{\deg(\omega_0') + \dots + \deg(\omega_n') - 3(n-m) + 2 \lvert \nu \rvert - m}{2} + k \\
                & \qquad = \dfrac{m_n}{2}+ k + m + \lvert \nu \rvert - n \\
                & \qquad \leq l - \frac{m + \lvert \nu \rvert}{2}+ m + \lvert \nu \rvert - n \\
                & \qquad = l + \dfrac{m + \lvert \nu \rvert}{2}- n \\
                & \qquad \leq l + m - n \leq l - 1.
		\end{align*}
		Therefore, the first part of Patodi's lemma implies that
		\begin{align*}
                & \lim_{t \to 0+}F_{x;t}(\omega_{0}, \ldots, \omega_{n}) \\
                & \quad = (-2)^{n} (2 \pi)^{-l}2^{-m_n/2}\lim_{t \to 0+}t^{\frac{m_n}{2} - l} \Str_x\left(\mathbf{c}_{x}(\omega_{0}') 
                \int_{t \sigma_n}\delta \Psi^{x;t; 1}_{\omega_1}(s_{1}) \cdots \delta \Psi^{x;t; 1}_{\omega_n}(s_{n}) \, \sum_{k = 0}^{l} \frac{T_x(\mathsf{A}^x(t))^{k}}{k!}\right),
		\end{align*}
	in case the latter limit exists. Recall that
		\begin{equation*}
			\Psi^{x;t; 1}_{\omega_j}(s_{j})=\int^{s_j}_{0}\pa^{x;t}(s)^{-1}V_{P(\omega_j)}(\mathsf{B}^{x;t}_{s})\pa^{x;t}(s) \, \Id s
		\end{equation*}
		and
		\begin{equation*}
			V_{P(\omega_j)}= \mathbf{c}(\Id^{\dagger} \omega_{j}') + \mathbf{c}(\omega_{j}'').
		\end{equation*}
		Let $l^{*} \coloneqq l - n - \frac{m_{n}}{2}$ and $0 \leq k \leq l^{*} - 1$. It follows that the Patodi order of every
		summand in
		\begin{equation*}
			\mathbf{c}_{x}(\omega_{0}') \int_{t \sigma_n}\delta \Psi^{x;t; 1}_{\omega_1}(s_{1}) \cdots \delta \Psi^{x;t; 1}_{\omega_n}
			(s_{n}) \, \frac{T_x(\mathsf{A}^x(t))^{k}}{k!}
		\end{equation*}
		is less or equal to
		\begin{equation*}
			\frac{m_{n} + 2n}{2}+ k \leq l - 1,
		\end{equation*}
		so that
		\begin{equation*}
			\ldots = (-2)^{n} (2 \pi)^{-l}2^{-m_n / 2}\lim_{t \to 0+}t^{\frac{m_n}{2} - l} \Str_x\left(\mathbf{c}
			_{x}(\omega_{0}') \int_{t \sigma_n}\delta \Psi^{x;t; 1}_{\omega_1}(s_{1}) \cdots \delta \Psi^{x;t; 1}_{\omega_n}(s_{n}
			) \, \frac{T_x(\mathsf{A}^x(t))^{l*}}{l^{*}!}\right),
		\end{equation*}
  in case the latter limit exists. Similarly to the above estimates of Patodi orders, we see that the only summand in
		\begin{equation*}
			\int_{t \sigma_n}\delta \Psi^{x;t; 1}_{\omega_1}(s_{1}) \cdots \delta \Psi^{x;t; 1}_{\omega_n}(s_{n})
		\end{equation*}
		contributing to the limit is the one where every integrand is of the form $\mathbf{c}(\omega_{j}'')$ for
		$1 \leq j \leq n$. Since we have an ordinary $n$-fold Riemann integral, it follows from by Lemma \ref{Limit parallel transport} that
		\begin{align*}
                \ldots & = (-2)^{n} (2 \pi)^{-l}2^{-\frac{m_n}{2}} \Str_x\left(\mathbf{c}_{x}(\omega_{0}') \lim_{t \to 0+}t^{-n}\int_{t \sigma_n}\delta 
                \Psi^{x;t; 1}_{\omega_1}(s_{1}) \cdots \delta \Psi^{x;t; 1}_{\omega_n}(s_{n}) \, \frac{T_x(\mathsf{A}^x(t))^{l^*}}{t^{l^*}l^{*}!}\right) \\
                & = \dfrac{(-1)^n 2^{-\frac{m_n}{2} + n}}{n!(2 \pi)^l} \Str_x\left(\mathbf{c}_{x}(\omega_{0}') \mathbf{c}_{x}(\omega_{1}'') \cdots \mathbf{c}_{x}(\omega_{n}'') \, \frac{T_x(\mathsf{D}^x)^{l^*}}{l^{*}!}\right).
            \end{align*}
            As above, one calculates
            \begin{align*}
                & \mathbb{E}[F_x(\omega_0, \dots, \omega_n)] \cdot \underline{\mu}_x \\
                & \quad = \dfrac{(-1)^n 2^{-\frac{m_n}{2} + n}}{n!(2 \pi)^l} \mathbb{E}\left[\Str_x\left(\mathbf{c}_{x}(\omega_{0}') \mathbf{c}_{x}(\omega_{1}'') \cdots \mathbf{c}_{x}(\omega_{n}'') \, 
                \frac{T_x(\mathsf{D}^x)^{l^*}}{l^{*}!}\right)\right]e^1_x \wedge\cdots\wedge e^d_x  \\
                & \quad = \dfrac{(-1)^n 2^{2n}}{n!(2\pi \sqrt{-1})^l}\big(\hat{A}(X) \wedge \omega_{0}' 
                \wedge \omega_{1}'' \wedge \cdots \wedge \omega_{n}''\big)^{[d]}_{x},
            \end{align*}
		which concludes the proof.
        \end{proof}

	\section{appendix}
	\label{appe}

	Let $X$ be a (for simplicity) closed and connected Riemannian manifold with $p(t,x,y)>0$ the heat kernel of the scalar
	Laplacian $-\Delta/2$ with respect to the volume measure $\mu$. Given $t>0$ we denote with $C([0,t],X)$ the Wiener
	space of continuous paths $\gamma\colon[0,t]\to X$. We give the latter the topology of uniform convergence. Let
	$\IFF^{X}_{t}$ denote the Borel-sigma algebra on $C([0,t],X)$.

	For every $x,y\in X$, $t>0$, the \emph{pinned Wiener measure $\IP^{x,y;t}$ from $x$ to $y$ with terminal time $t$} is the
	unique probability measure on $\IFF^{X}_{t}$, such that for all $0<t_{1}<\dots t_{n}<t$, and all Borel sets $U_{1},\dots, U_{n}\subset X$,
	\begin{align*}
            & p(t,x,y)\cdot\IP^{x,y;t}\{\gamma\colon \gamma(t_{1})\in U_{1},\dots,\gamma(t_{n})\in U_{n}\} \\
            & \qquad = \int_{U_1}\cdots\int_{U_n}p(t-t_{n},y,x_{n})p(t_{n}-t_{n-1},x_{n},x_{n-1}) \times \\
            & ~~~~~~~~~~~~~~~~~~ \times \cdots p(t_{2}-t_{1},y,x_{n})p(t_{1},x_{1},x)\, \Id\mu(x_{1})\cdots \Id\mu(x_{n}).
	\end{align*}
	A detailed proof of the existence of this measure can be found in \cite{bar}. One has the normalization
	\begin{equation*}
            \IP^{x,y;t}\{\gamma\colon\gamma(0)=x,\gamma(t)=y\}=1.
	\end{equation*}

	A process $\mathsf{B}^{x,y;t}_{s}\in X$, $s\in [0,t]$, which is defined on a probability space with filtration $(\IFF_{s}
	)_{s\in [0,t]}$, is called an \emph{adapted Brownian bridge from $x$ to $y$ with terminal time $t$}, if
	\begin{itemize}
		\item the law of $\mathsf{B}^{x,y;t}$ is equal to $\IP^{x,y;t}$,

		\item $\mathsf{B}^{x,y;t}$ is adapted to $(\IFF_{s})_{s\in [0,t]}$,

		\item $\mathsf{B}^{x,y;t}$ satisfies the following time-inhomogeneous Markov property: for all $s\in [0,t)$, all $\IFF
			_{s}$-measurable random variables $\Phi\colon\Omega \to [0,\infty]$, and for all random variables
			$\Psi\colon C([0,t],X)\to [0,\infty]$ one has
			\begin{equation*}
				\mathbb{E}\left[ \Phi\cdot\left( \Psi \circ \mathsf{B}^{x,y;t}_{\min(s+\bullet,t)}\right)\right]=\mathbb{E}\left[
				\Phi\cdot\int_{C([0,t],X)}\Psi\Big(\gamma(\min(\bullet,t-s)\Big) \, \Id\IP^{\mathsf{B}^{x,y;t}_s,y;t-s}(\gamma)\right].
			\end{equation*}
	\end{itemize}

	A fundamental result (which can be traced back at least to \cite{bismut} in the compact case; see also \cite{batu2})
	states that any adapted Brownian bridge as above is a continuous semimartingale on the whole time interval $[0,t]$, meaning
	by definition that for all smooth functions $f\colon X\to\IR$ the real-valued process $f(\mathsf{B}^{x,y;t})$ is a continuous
	semimartingale.

\end{document}